\documentclass{siamart1116}
\usepackage{graphicx, epsfig, psfrag}
\usepackage{amsfonts}
\usepackage{amssymb}
\usepackage{amsmath}
\usepackage{cases}
\usepackage{mathrsfs}
\usepackage{subfig}
\usepackage{bm}
\usepackage{color}
\usepackage{float}
\usepackage{hyperref}
\usepackage{mathrsfs}
\usepackage{longtable}
\usepackage[top=3cm,bottom=3cm,left=3cm,right=3cm]{geometry}
\usepackage{epstopdf}
\usepackage{mathabx}
\usepackage{stmaryrd}
\usepackage{ulem}
\usepackage{scalerel}
\usepackage{enumitem}
\usepackage{xcolor}
\usepackage{lipsum}
\usepackage{amsfonts}
\usepackage{graphicx}
\usepackage{epstopdf}
\usepackage{algorithmic}

\ifpdf
  \DeclareGraphicsExtensions{.eps,.pdf,.png,.jpg}
\else
  \DeclareGraphicsExtensions{.eps}
\fi

\numberwithin{theorem}{section}

\newcommand{\TheTitle}{A Posteriori Error Estimation and Adaptive Algorithm for Atomistic/Continuum Coupling in 2D} 
\newcommand{\TheAuthors}{H. Wang, M.Liao, P.Lin and L. Zhang}
\newcommand{\TheSubtitle}{A Posteriori Estimation A/C Coupling 2D}
\headers{\TheSubtitle}{\TheAuthors}

\title{{\TheTitle}\thanks{Submitted to the editors on \today.
\funding{HW was partially supported by NSFC grant 11501389, 11471214 and Sichuan University Starting Up Research Funding No. 2082204194117. PL and ML were partially supported by NSFC grant 91430106 and Fundamental Research Funds for Central Universities Nos. 06108038 and FRF-BR-13-023. PL and HW were partially supported by EMS RSF grant. LZ was partially supported by NSFC grant 11471214, 11571314 and the One Thousand Plan of China for young scientists.}}}

\author{
  Hao Wang\thanks{College of Mathematics, Sichuan University, No.24 South Section One, Yihuan Road, Chengdu, 610065, China
    (\email{wangh@scu.edu.cn}).}
  \and
  Mingjie Liao\thanks{Department of Applied Mathematics and Mechanics,
University of Science and Technology Beijing,
No. 30 Xueyuan Road, Haidian District,
Beijing 100083 
(\email{mliao@xs.ustb.edu.cn}).}
  \and
  Ping Lin\thanks{Department of Mathematics,
University of Dundee,
Dundee, DD1 4HN, Scotland,
United Kingdom
(\email{plin@maths.dundee.ac.uk}).}
  \and
  Lei Zhang\thanks{Institute of Natural Sciences, School of Mathematical Sciences,
  and Ministry of Education Key
  Laboratory of Scientific and Engineering Computing (MOE-LSC),
  Shanghai Jiao Tong University, 800 Dongchuan Road, Shanghai
  200240, China
(\email{lzhang2012@sjtu.edu.cn}).}
}

\usepackage{amsopn}

\ifpdf
\hypersetup{
  pdftitle={\TheTitle},
  pdfauthor={\TheAuthors}
}
\fi

\renewcommand{\cases}[1]{\left\{ \begin{array}{rl} #1 \end{array} \right.}

\newcommand{\mymat}[1]{\left[ \begin{matrix} #1 \end{matrix} \right]}

\def\XXint#1#2#3{{\setbox0=\hbox{$#1{#2#3}{\int}$ }
\vcenter{\hbox{$#2#3$ }}\kern-.6\wd0}}
\def\b{\big}
\def\B{\Big}
\def\bg{\bigg}

\def\sep{\,|\,}
\def\bsep{\,\b|\,}

\DeclareMathOperator*{\argmin}{argmin}

\def\supp{{\rm supp}}

\def\R{\mathbb{R}}

\def\Z{\mathbb{Z}}

\def\bbV{\mathbb{V}}

\def\dx{\,{\rm d}x}

\def\ds{\,{\rm d}s}

\def\pp{\partial}

\def\db{\,{\rm db}}

\def\<{\langle}
\def\>{\rangle}

\def\mA{{\sf A}}
\def\mB{{\sf B}}

\def\mF{{\sf F}}
\def\mG{{\sf G}}
\def\mH{{\sf H}}
\def\mI{{\sf I}}
\def\mJ{{\sf J}}
\def\mP{{\sf P}}

\def\mR{{\sf R}}

\def\bbB{\mathbb{B}}

\def\Hs{\mathcal{H}}

\newcommand{\Dc}[1]{\D_{#1}}
\def\D{\nabla}
\def\del{\delta}
\def\ddel{\delta^2}

\def\a{{\rm a}}

\def\c{{\rm c}}

\def\a{{\rm a}}
\def\c{{\rm c}}

\def\i{{\rm i}}

\def\supp{{\rm supp}}

\def\L{\Lambda}
\def\Cs{\mathcal{C}}
\def\Fs{\mathcal{F}}

\def\Us{\mathscr{U}}

\def\Ush{\dot{\Us}^{1,2}}

\def\Usc{\Us^c}

\def\Bs{\mathcal{B}}

\def\E{\mathscr{E}}
\def\Ea{\E^\a}
\def\Eb{\E^{\rm b}}

\def\Ec{\E^\c}
\def\Eh{\E^\h}

\def\L{\Lambda}
\def\La{\L^{\a}}
\def\Li{\L^{\i}}

\def\Nhd{\mathcal{N}}
\def\rcut{r_{\rm cut}}
\def\Rg{\mathscr{R}}

\def\vsig{\varsigma}
\def\T{\mathcal{T}}
\def\Tp{T^+}
\def\Tm{T^-}
\def\np{\nu^+}
\def\nm{\nu^-}
\def\Th{\mathcal{T}_h}

\def\Fh{\mathscr{F}_h}
\def\UsT{\Us_h}

\def\Ia{I_\a}

\def\vor{\rm vor}
\def\s{\sigma}
\def\sh{\sigma^h}
\def\sa{\sigma^\a}
\def\sigc{\sigma^\c}
\def\Oma{\Omega^\a}
\def\PO{{\rm P}_0}
\newsiamremark{remark}{Remark} 
\newsiamremark{assumption}{Assumption} 
\definecolor{lzcol}{rgb}{0.7, 0, 0}
\definecolor{hwcol}{rgb}{0, 0, 0.9}
\definecolor{mlcol}{rgb}{0, 0.7, 0}
\definecolor{todocol}{rgb}{0.0, 0.4, 0.0}

\begin{document}

\maketitle

% REQUIRED
\begin{abstract}
Atomistic/continuum coupling methods aim to achieve optimal balance between accuracy and efficiency. Adaptivity is the key for the efficient implementation of such methods. In this paper, we carry out a rigorous a posteriori analysis of the residual, the stability constant, and the error bound, for a consistent atomistic/continuum coupling method in 2D.  We design and implement the corresponding adaptive mesh refinement algorithm, and the convergence rate with respect to degrees of freedom is optimal compare with a priori error estimates. 
\end{abstract}

% REQUIRED

% REQUIRED

\maketitle

%===========================================================================

\section{Introduction}
\label{sec:introduction}
Atomistic/continuum (a/c) coupling methods are a class of computational multiscale methods that aim to combine the accuracy of the atomistic model and the efficiency of the continuum model for crystalline solids with defects \cite{Ortiz:1995a, Shenoy:1999a, Gumbsch:1989}. Namely, the atomistic model can be applied in a small neighborhood of the localized defects such as vacancies, dislocations, and cracks, while the continuum model (e.g., Cauchy-Born rule) can be employed away from the defect cores where elastic deformation occurs. The construction and analysis of different a/c coupling methods have attracted considerable attention in the research community in recent years \cite{LinP:2006a, OrtnerShapeev:2011, Luskin:clusterqc, MiLu:2011}. We refer the readers to \cite{Miller:2008, LuOr:acta} for a review of such methods. 

The goal of the mathematical analysis for a/c coupling methods is to find the optimal relation of accuracy vs. degrees of freedom. The a priori analysis has been carried out for several typical a/c coupling methods, for example the QNL (quasi-nonlocal quasicontinuum) method \cite{emingyang, OrtnerWang:2011}, the BQCE (blended energy-based quasi-continuum) method \cite{LiOrShVK:2014}, the BQCF (blended force-based quasi-continuum) method \cite{MiLu:2011, LiOrShVK:2014}, the GRAC (geometric reconstruction based atomistic/continuum coupling) method \cite{PRE-ac.2dcorners} and the BGFC (atomistic/continuum blending with ghost force correction) method \cite{OrZh:2016}. 

In contrast,  although adaptivity is the key for the efficient implementation of a/c coupling methods, only few research articles are concerned with the a posteriori error control of these methods. The goal-oriented approach has been utilised in \cite{prud06} by Prudhomme et al. to provide a posteriori error control for a three dimensional nanoindentation problem with the quantity of interest being the force acting on the indenter. The error estimator is a modification of the rigorously derived residual functional, and its effectiveness is only validated numerically. Arndt and Luskin \cite{arndtluskin07b, arndtluskin07c} analyze the goal-oriented approach for a one dimensional Frenkel-Kontorova model, where the a posteriori error estimators are used to optimize the choice of the atomistic region as well as the finite element mesh in the continuum region. All these work employ the original energy-based quasicontinuum method as the underlying model which is later shown to be inconsistent and suffers from the so-called "ghost force"  \cite{Shenoy:1999a, Dobson:2008c, LinShapeev:2009, emingyang, Makridakis2012}. Recently, Kochmann et al. \cite{Kochman:2016} proposed an adaptivity strategy for the so-call "fully-nonlocal quasi-continuum" method which apply a discrete model in the entire computational domain without coupling of different models. This approach aims to minimize the ghost force rather than eliminate it as in the consistent a/c coupling method.

The residual based a posteriori error bounds for a/c coupling schemes are first derived in \cite{OrtnerSuli:2008a, Ortner:qnl.1d} by Ortner et al. in 1D. A recent advance in this direction \cite{OrtnerWang:2014} is the a posteriori error analysis of a consistent energy-based coupling method developed in \cite{Shapeev:2010a, Shapeev2012}, where the a posteriori error estimators are proposed both in the energy norm and in energy itself.  For complex lattice, a posteriori error analysis for the QC method in 1D has been carried out in \cite{Abdulle:2013}.

Despite all those developments, the rigorous mathematical justification of a posteriori error estimates beyond 1D is still missing. In this paper, we present a rigorous a posteriori error estimate for a consistent energy-based a/c method in two dimension, which is of physical significance and has not been considered so far to the best knowledge of the authors. We use the residual-based approach \cite{Verfurth:1996a} to establish the estimate in negative Sobolev norms following \cite{OrtnerWang:2014}. Two features distinguish our problem from the classic residual-based estimate for finite element approximation of the elliptic equations. The first one is the existence of the modeling error which is in origin different from the applications of quadrature rules. The second one is that the mesh may not be further refined when it almost coincides with the reference lattice, therefore a model adaptation should be imposed.  The analysis and algorithm rely on the so-called divergence free tensor field, which characterizes the essential difference of 2D results compared with 1D results in \cite{Ortner:qnl.1d, OrtnerWang:2014} where the analysis can be carried out by explicit calculations.

Similar to the a priori analysis of GRAC in \cite{PRE-ac.2dcorners}, we constrain ourselves to the case of nearest-neighbor interactions. Although the analysis can be extended to finite range interactions and to other a/c coupling methods, we decide not to include these so that the main ideas and steps are clearly presented without the distraction from the unnecessary complexity of the presentation. Instead, we will make further remarks on this point again in \S~\ref{sec:conclusion}. 

The paper is organized as follows. In \S~\ref{sec:formulation} we set up the atomistic, continuum and coupling models for point defects. In \S~\ref{sec:error} we present the main results: the residual estimate, stability bound, and rigorous a posteriori error estimates for the coupling scheme. We formulate the corresponding adaptive algorithm and demonstrate numerical results in \S~\ref{sec:numerics}. We draw conclusions and make suggestions for future research in \S~\ref{sec:conclusion}. Some auxiliary results are given in \S~Appendix~\ref{sec:appendix:extension}.

\section{Formulation}
\label{sec:formulation}

We first give a brief review of a model for crystal defects in an infinite lattice in the spirit of \cite{2013-defects} in \S~\ref{sec:formulation:atm} and the Cauchy-Born continuum model in \S~\ref{sec:formulation:cb}. We then present a generic form of a/c coupling schemes in \S~\ref{sec:formulation:ac}. We will introduce the consistent scheme GRAC specifically in \S~\ref{sec:formulation:grac}. 

\subsection{Atomistic model}
\label{sec:formulation:atm}
\def\Rdef{R^{\rm def}}
\def\Rg{\mathcal{R}}
\def\Rgnn{\mathcal{N}}
\def\rcut{r_{\rm cut}}
\def\Lhom{\L^{\rm hom}}
\def\Ddef{D^{\rm def}}
\def\Ldef{\L^{\rm def}}
\subsubsection{Atomistic lattice and defects}
\label{sec:formulation:atm:lattice}

Given $d \in \{2, 3\}$, $\mA \in \R^{d \times d}$ non-singular,  $\Lhom := \mA \Z^d$ is the \textit{homogeneous reference lattice} which represents a perfect \textit{single lattice crystal} formed by identical atoms and possessing no defects. $\L\subset \R^d$ is the \textit{reference lattice} with some local \textit{defects}. The mismatch between $\L$ and $\Lhom$ represents possible defects $\Ldef$, which are contained in some localized \textit{defect cores} $\Ddef$ such that the atoms in $\L\setminus \Ddef$ do not interact with defects $\Ldef$ (see \S~\ref{sec:formulation:atm:function} and \S~\ref{sec:formulation:atm:potential} regarding interaction neighbourhood). For example, $\Ldef = \{x\}$ for a crystal with a single point defect at $x$, and one can choose a proper radius $\Rdef > 0$ such that $\Ddef = B_{x, \Rdef}$, where $B_{x, R} := \{ z \in \R^d \sep |z-x| \leq R\}$.  For different types of point defects, we have
\begin{itemize}
\item $\L\subset \Lhom$ for a vacancy at $x\in \Lhom$; 
\item $\L \supset \Lhom$ for an interstitial at $x\in \L$ but $x\notin\Lhom$; 
\item $\L = \Lhom$ for an impurity at $x\in\Lhom$, the difference of the impurity atom with other atoms can be characterized by the inhomogeneity of interaction potentials (see \S~\ref{sec:formulation:atm:potential}). 
\end{itemize}
This characterization of localized defects can be straightforwardly generalized to multiple point defects and micro-cracks, for example, see the setup of the model problem in \S~\ref{sec:numerics:problem}. Straight screw dislocations can be enforced through the appropriate choice of boundary conditions \cite{2013-defects}.

%For all point defects, $d = m \in \{2, 3\}$, and one takes $\Phi_\ell$
%to be homogeneous outside some core radius, that is, $\Phi_\ell(y) =
%\Phi(y-y(\ell))$ for some site potential $\Phi$. A suitable boundary
%condition in this case is $\yd(x) = \mB x$, for some ``macroscopic
%strain'' $\mB \in \R^{d \times d}$. This satisfies all conditions of
%the theory.

% These models are, strictly speaking, for dimension $d = 3$, but
% ``toy models'' can also be formulated for $d = 2$.

\subsubsection{Lattice function and lattice function space}
\label{sec:formulation:atm:function}

Given $d\in\{2,3\}$, $m\in\{1,2,3\}$, denote the set of vector-valued \textit{lattice functions} by 
\begin{displaymath}
\Us := \{v: \L\to \mathbb{R}^m \}.
\end{displaymath}

A \textit{deformed configuration} is a \textit{lattice function} $y \in \Us$. Let $x$ be the identity map, the \textit{displacement} $u\in \Us$ is defined by $u(\ell) = y(\ell)-x(\ell) = y(\ell) - \ell $ for any $\ell\in\L$.

For each $\ell\in \L$, we prescribe an \textit{interaction neighbourhood}  $\Nhd_{\ell} := \{ \ell' \in \L \sep 0<|\ell'-\ell| \leq \rcut \}$ with some cut-off radius $\rcut$. The \textit{interaction range} $\Rg_\ell := \{\ell'-\ell \sep \ell'\in \Nhd_\ell\}$ is defined as the union of lattice vectors defined by the finite difference of lattice points in $\Nhd_{\ell}$ and $\ell$. 

To measure the error for lattice functions we need to introduce function norms and function spaces on the lattice. Define the ``finite difference stencil'' $Dv(\ell):= \{D_\rho v(\ell)\}_{\rho \in \Rg_\ell} :=\{v(\ell+\rho)-v(\ell)\}_{\rho \in \Rg_\ell}$. Higher-order finite differences, e.g., $D_\rho D_\vsig v$ and $D^2 v$ can be defined in a canonical way. A \textit{lattice function norm} can hence be defined using those notations. For $v \in\Us$, let the lattice energy-norm (a discrete $H^1$-semi-norm) be
\begin{equation}
  \| Dv \|_{\ell^2} := \bg( \sum_{\ell \in \L}
  \sum_{\rho \in \Rg_\ell} |D_\rho v(\ell)|^2
  \bg)^{1/2}.
\end{equation}
The associated \textit{lattice function space} is defined by
\begin{align*}
  \Ush &:= \b\{ u : \L \to \R^m \bsep \| Du \|_{\ell^2} < +\infty \b\}.
\end{align*}

We choose
\begin{equation}
\Bs: = \{(\ell, \ell+\rho): \ell \in \L, \rho \in \Rg_\ell\}
\label{eq:bonds}
\end{equation}
to be the collection of all the nearest neighbour bonds in the reference lattice, and for $b=(\ell, \ell+\rho)\in \Bs$, denote $\rho_b = \rho$. Then the energy norm can be reformulated as 
\begin{equation}
  \| Dv \|_{\ell^2} := \bg( \sum_{b=(\ell, \ell+\rho) \in \Bs} |D_{\rho} v(\ell)|^2
  \bg)^{1/2}.
\end{equation}

The homogeneous lattice $\Lhom = \mA\mathbb{Z}^d$ naturally induces a simplicial micro-triangulation $\T$. In 2D, $\T^{\a} = \{\mA\xi + \hat{T}, \mA\xi-\hat{T}|\xi\in \mathbb{Z}^2\}$, where $\hat{T} = {\rm conv}\{0, e_1, e_2\}$. Let $\bar{\zeta}\in W^{1, \infty}(\Lhom; \R)$ be the P1 nodal basis function associated with the origin; namely, $\bar{\zeta}$ is piecewise linear with respect to $\T^\a$, and $\bar{\zeta}(0) = 1$ and $\bar{\zeta}(\xi)=0$ for $\xi\neq 0$ and $\xi\in \Lhom$. The nodal interpolant of $v\in \Us$ can be written as 
\begin{displaymath}
	\bar{v}(x):=\sum_{\xi\in\Z^d}v(\xi)\bar{\zeta}(x-\xi).
\end{displaymath}

\def\Use{\Us^{1,2}}

We can introduce the discrete homogeneous Sobolev spaces
\begin{displaymath}
	\Use :=\{u\in \Us|\nabla \bar{u}\in L^2\},
\end{displaymath}
with semi-norm $\|\nabla \bar{u}\|_{L^2}$.  It is known from \cite{OrShSu:2012} that $\Ush$ and $\Use$ are equivalent.

\subsubsection{Interaction potential}
\label{sec:formulation:atm:potential}
%

% Multibody potential

For each $\ell \in \L$, let $V_\ell(y)$ denote the \textit{site energy}
associated with the lattice site $\ell \in \L$, and we assume that $V_\ell(y)\in C^k((\R^d)^{\Rg_\ell}), k \geq 2$. In this paper, we consider the general
multibody interaction potential of the \textit{generic pair functional form} \cite{TadmorMiller:2012}. Namely, the potential is a function of the distances between atoms within interaction range and with no angular dependence. Accordingly, we have the following equivalent forms of interaction potentials of generic pair functional form,
\begin{equation}
\label{eq:pairfunctional}
	V_{\ell}(y) = \widehat{V}_{\ell}(\{D_\rho y(\ell) \}_{\rho\in \Rg_\ell}) = \widetilde{V}_{\ell}(\{|D_\rho y(\ell) |\}_{\rho\in \Rg_\ell})
\end{equation}

\begin{remark}
\label{rem:pairfuncational}
For convenience, with a slight abuse of notation, we will use $V_{\ell}(D_\rho y)$,  $V_{\ell}(|D_\rho y|)$ instead of $\widehat{V}_{\ell}(\{D_\rho y(\ell) \}_{\rho\in \Rg_\ell})$,  $\widetilde{V}_{\ell}(\{|D_\rho y (\ell) |\}_{\rho\in \Rg_\ell})$ when there is no confusion in the context.
\end{remark}

We assume that $V_\ell$ is \textit{homogeneous} outside the defect region $\Ddef$, namely, $V_\ell = V$ and $\Rg_\ell = \Rg$ for $\ell \in \Lambda \setminus \Ddef$. $V$ and $\Rg$ have the following point symmetry: $\Rg = -\Rg$, and $V(\{-g_{-\rho}\}_{\rho\in\Rg}) = V(g)$.

\begin{remark}
\label{rmk:cfuy}
Notice that both displacement $u$ and deformation $y$ are discrete functions belonging to $\Us$, however $u\in\Ush$ while $y\notin \Ush$. We define the interaction potential $V$ through $y$ for the convenience of stability analysis, the consistency results are the same either with $u$ or with $y$.
\end{remark}

A great number of practical potentials are in the form \eqref{eq:pairfunctional}, including the widely used embedded atom model (EAM) \cite{Daw:1984a} and Finnis-Sinclair model \cite{FS:1984}. For example, assuming a finite interaction neighborhood $\Nhd_{\ell}$ and an interaction range $\Rg_{\ell}$ for $\ell\in\L$, EAM potential reads
\begin{align}
  V_{\ell}(y) := & \sum_{\ell' \in \Nhd_{\ell}} \phi(|y(\ell)-y(\ell')|) + F\B(
  {\textstyle \sum_{\ell' \in \Nhd_{\ell}} \psi(|y(\ell)-y(\ell')|)} \B),\nonumber\\
    = &\sum_{\rho \in \Rg_{\ell}} \phi\b(|D_\rho y(\ell)|\b) + F\B(
  {\textstyle \sum_{\rho \in \Rg_{\ell}}} \psi\b( |D_\rho y(\ell)|\b) \B).   
    \label{eq:eam_potential}
\end{align}
for a pair potential $\phi$, an electron density function $\psi$ and
an embedding function $F$. 

%To describe, e.g., impurity defects, we allow $\phi, G, \rho$ to be
%species-dependent, i.e., $\phi = \phi_{ab}, G = G_a, \psi =
%\psi_{ab}$.

%\def\Ldef{\L_{\rm def}}
%\def\Ldefhom{\L^{\rm hom}_{\rm def}}
%\def\Lren{L^{\rm ren}} 
%%
%We introduce the concept of \textit{homogeneity} for site potentials.  For $\ell \in \Lhom$ we say that the site potentials are \textit{globally homogeneous} if $V_{\ell}(y) = V_{\ell'}(z)$ for all $\ell, \ell' \in\Lhom$. In this case, it is easy to see (summation by parts, or point-symmetry of the lattice) that
%\begin{displaymath}
%  \sum_{a \in \Lhom} \b\< \del V_{\ell}(x), u \b\> = 0 \qquad \forall
%  u : \Lhom \to \R^d, {\rm supp}(u) \text{~compact}.
%\end{displaymath}
%
%For general $\L$, we assume that the site potentials are homogeneous
%outside the defect cores $\Ddef$. Let $V_{\ell}^{\rm
%  hom}, \ell \in \Lhom$ be a globally homogeneous site potential so that
%$\< \del V_{\ell}(x), u \> = \< \del V_{\ell}^{\rm hom}, u \>$ for all $\ell
%\in \L, |\ell| > \Rdef+\rcut$ and for all $u : \L \cup \Lhom \to \R^d$
%with compact support. Homogeneity of the site potential outside the defect cores $\Ddef$ entails
%simply that only one atomic species occurs there. 
%
% Well-posedness

The energy of an infinite configuration is typically ill-defined. However, if we redefine the potential $V_\ell(y)$ as the difference $V_\ell(y) - V_\ell(\ell)$, which is equivalent to assuming $V_\ell(\ell) =0$,  the energy functional
\begin{equation}
\label{eqn:Ea}
  \Ea(y) = \sum_{\ell \in \L} V_\ell(y) 
\end{equation}
is a meaningful object. Given the point symmetry and smoothness assumptions for the site potentials $V_\ell$, $\Ea(y)$ is well-defined for $y-y^{B}  \in \Use$, where $y^B(x)=Bx$. Furthermore, if $V_\ell(y)$ is $C^k$ in its variables,  $\Ea$ is $k$ times Fr\'{e}chet differentiable. In particular, we define $M$ as the Lipschitz constant of $\ddel\Ea$, by \cite[Lemma 2.1]{2013-defects}. 

Under the above conditions, the goal of the atomistic problem is to find a \textit{strongly stable}
equilibrium $y$, such that, given a macroscopic applied strain $B\in\R^{d\times d}$, we aim to compute
\begin{equation}
  \label{eq:min}
  y \in \arg\min \b\{ \Ea(y) \bsep y-y^B\in \Use \b\}.
\end{equation}
$y$ is \textit{strongly stable} if there exists $c_0 > 0$ such that
\begin{displaymath}
\< \ddel \Ea(y) v, v \> \geq c_0 \| \nabla v \|_{L^2}^2, \quad \forall v \in \Us^{1,2}.
\end{displaymath}.

%\begin{remark}
%  The condition $\del E(\yd;\yd) \in \Ushd$ encodes that $\yd$ is an
%  approximate equilibrium in the far-field. It is easy to see that,
%  without such a condition, $\inf_{\| u \|_{\Ush} \leq \epsilon}
%  E(y+u;\yd) = -\infty$ for all $\epsilon > 0$ and for all $y \in
%  \yd+\Ush$. That is, there can exist no local energy minimizers of
%  \eqref{eq:min}.
%\end{remark}
%
%Under suitable conditions on the site potentials $V_{\ell}, a \in \L$
%(most crucially, regularity and homogeneity outside the defect;
%cf. \S~\ref{sec:scheme:hom}), it is shown in
%\cite{2013-defects} that the energy-difference functional
%\begin{equation}
%  \label{eq:scheme:Ea}
%  \Ea(u) := \sum_{a \in \L} \Phi'_a(u), \qquad \text{where } \Phi'_a(u)
%  := V_{\ell}(x+u)-V_{\ell}(x),
%\end{equation}
%is well-defined for all relative \textit{displacements} $u \in \Us$,
%where $\Us$ is given by
%\begin{align*}
%  \Us &:= \b\{ v : \L \to \R^d \bsep |v|_{\Us} < +\infty \b\}, \quad
%  \text{where} \\
%  |v|_{\Us} &:= \B(\sum_{a \in \L} \sum_{b \in \Rgnn_a} |v(b)-v(a)|^2 \B)^{1/2}.
%\end{align*}
%where $\Rgnn_a\subset \L\setminus\{a\}$ is a set of
%``nearest-neighbour'' directions satisfying ${\rm span}(\Rgnn_a) =
%\R^d$ and $\sup_{a \in \L} \# \Rgnn_a < \infty$ for $a \in \L$.  

%\def\ua{u^\a}
%
%\lz{
%Leave the external force as a remark???
%
%Under external force $f$, the atomistic problem is to compute
%\begin{equation}
%  \label{eq:atm}
%  \ua \in \arg\min \b\{ \Ea(v) - \F(v) \bsep v \in \Us \b\}.
%\end{equation}
%where $\F(v):=\sum_{a\in\Lambda} f v$.
%This problem is analysed in considerable detail in
%\cite{2013-defects}.
%}

It is proven in \cite[Theorem 2.3 ]{2013-defects} that, if the homogeneous lattice is stable and $y\in \Us$ is a critical point of $\Ea$ such that $u=y-y^B\in\Use$, then $D^j u $ exhibit the following \textit{generic decay}, $j = 0, 1, \dots$,
\begin{equation}
  \label{eq:ptdef:regularity}
  \b|D^j u(\ell)\b| \lesssim |\ell|^{1-d-j}, \quad \text{and}
  \quad
  \b|u(\ell) - u_\infty \b| \lesssim |\ell|^{-d+1}.
\end{equation}
where $u_\infty:=\lim_{|\ell|\to \infty} u(\ell)$. 

\subsection{Continuum model}
\label{sec:formulation:cb}
To formulate atomistic to continuum coupling schemes, we need a
continuum model which is compatible with \eqref{eqn:Ea} and defined through a
strain energy density function $W : \R^{d \times d} \to \R$. Let $V$ 
be the homogeneous site potential on $\Lhom$. A typical choice
in the multi-scale context is the Cauchy--Born continuum model \cite{E:2007a, OrtnerTheil2012}, 
the energy density $W$ is defined by
\begin{displaymath}
  W(\mF) := \det \mA^{-1} V(\mF x).
\end{displaymath}

\def\Th{\mathcal{T}_h}
\def\Nh{\mathcal{N}_h}
\def\Ush{\Us_h}
\def\Ra{R^\a}
\def\Rb{R^{\rm b}}
\def\Rc{R^\c}
\def\Eb{\E^{\rm b}}
\def\dof{{\rm DOF}}
\def\Omh{\Omega_h}
\def\Thr{{\T_{h,R}}}
\def\vor{\rm vor}
\def\Uhr{\Us_{h,R}}

\subsection{A/C coupling}
\label{sec:formulation:ac}
We give a generic formulation of the a/c coupling method and employ concepts and notation from various earlier works, such as \cite{Ortiz:1995a,Shenoy:1999a,Shimokawa:2004,2012-optbqce, COLZ2013}, and we adapt the formulation to the settings in this paper. 

First, the computational domain $\Omega_R \subset \R^d$ is a simply connected, polygonal and closed set, such that $B_{0, R} \subset \Omega_R \subset B_{0, c_0 R}$ for some $c_0 >0$. Let $R$ be the radius of $\Omega_R$ We have the following decomposition $\Omega_R = \Omega_R^\a \bigcup \Omega_R^\c$, where the atomistic region $\Omega^\a_R$ is again simply
connected and polygonal, and contains the defect core: $\Ddef \subset \Omega^\a_R$. Let $R_\a$ be the radius of $\Omega^\a_R$. Let $\T^{c}_{h,R}$ be a shape-regular simplicial partition (triangles for $d = 2$ or tetrahedra for $d = 3$) of the continuum region $\Omega^\c_R$.  

\begin{figure}[htb]
\begin{center}
	\includegraphics[scale=0.4]{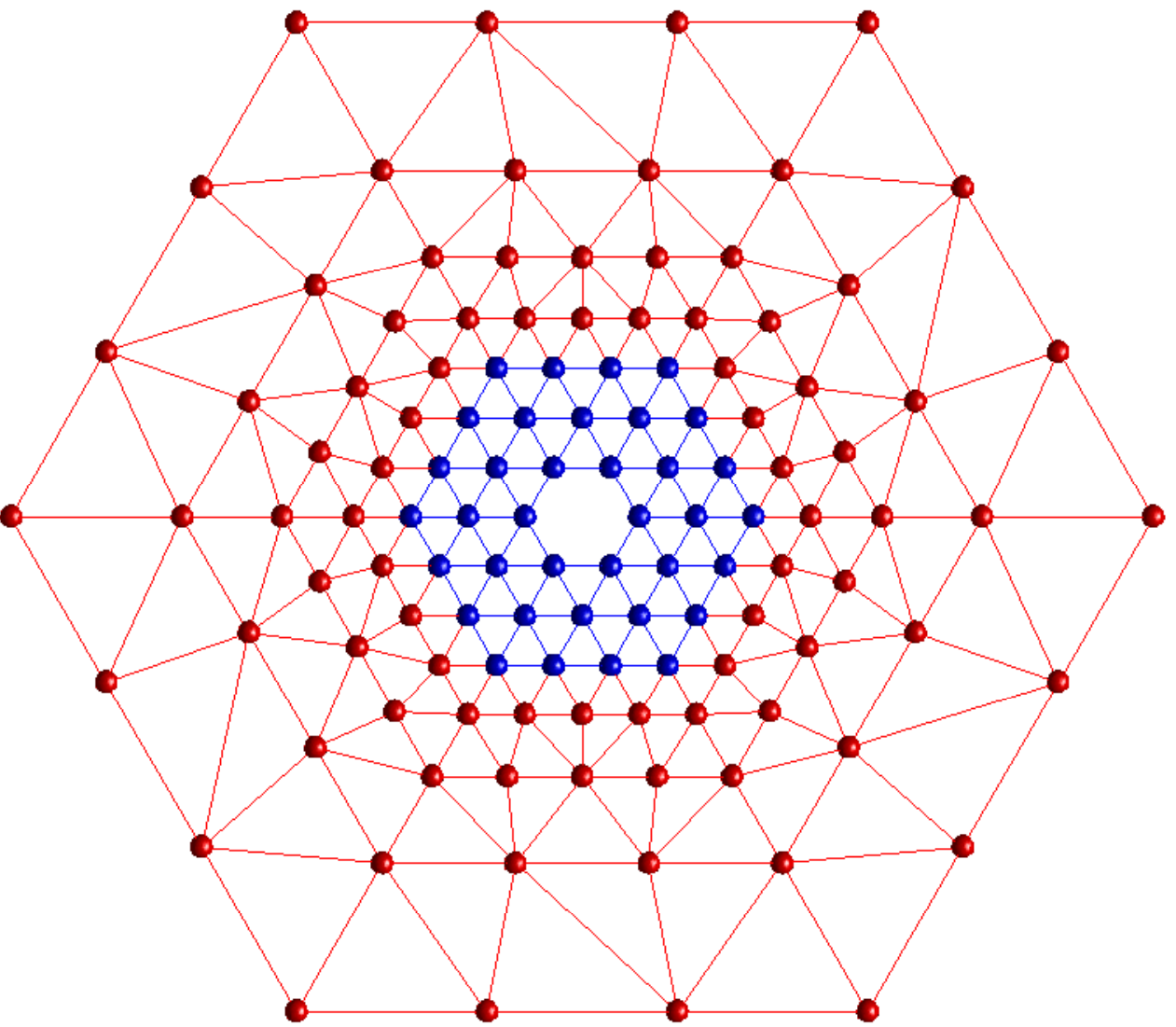}
	\caption{Illustration of computational mesh. The computational domain is $\Omega_R$, and the corresponding triangulation is $\T_{h,R}$. Blue nodes in $\Omega_R$ are atoms in $\L^{\a,\i}$. For nearest neighbour interaction, $\L^\i$ is the set of outmost layer of blue atoms. Red nodes in $\T_{h,R}$ are continuum degrees of freedom.  $\Omega_R^\a$ is the domain induced by the blue nodes, and $\T^\a_{h,R}$ is the corresponding triangulation. $\Omega_R^c$ and $\T^\c_{h,R}$ are the respective complements of $\Omega_R^\a$ and $\T_{h,R}$.}
	\label{figs:plotMesh}
\end{center}
\end{figure}

Next, we decompose the set of atoms $\L^{\a,\i} := \L \bigcap \Omega^\a_R = \L^\a \bigcup \L^\i$ into a core atomistic set $\L^\a$ and an interface set $\L^\i$ (typically a few ``layers'' of atoms surrounding $\L^\a$) such that $\L\bigcap\Ddef \subset \L^\a$. Let $\T^\a_{h,R}$ be the canonical triangulation induced by $\L^{\a,\i}$, which may contain "holes" due to the existence of defects, and $\Thr = \T^\c_{h,R} \bigcup \T^\a_{h,R}$. Sometimes, it is also convenient to define $\T^\i_{h,R}:=\{T\in\Thr: \L^{\i}\bigcap T \neq \emptyset\}$. Please see Figure \ref{figs:plotMesh} for an illustration of the computational mesh.

% fig illustration of mesh.

Let $\Omega_{h, R} = \bigcup_{T\in \Thr} T$. Notice that $\Omega_{h, R}$ can be multiple-connected, and $\Omega_R \setminus \Omega_{h, R}$ characterizes possible defects. The space of \textit{coarse-grained} displacements is,
\begin{align*}
  \Us_{h,R} := \b\{ u_h : \Omega_{h, R} \to \R^m \bsep ~&
  \text{ $u_h$ is continuous and p.w. affine w.r.t. $\T_{h,R}$, } \\[-1mm]
  & \text{ $u_h = 0$ on $\partial \Omega_R$ } \b\}.
\end{align*}

We may drop the subscript $R$ in the above definitions, for example, use $\Th$ instead of $\Thr$ if there is no confusion. Let $\Nh$ be the set of nodes in $\Th$, and $\mathcal{F}_h$ be the set of edges in $\Th$.

Denote $\vor(\ell)$ as the voronoi cell associated with atom $\ell$, the volume of this cell denoted as $|\vor(\ell)|$ equals the volume of the unit cell in $\Lhom$, i.e. $\vor(\ell) = \det (\mA)$. For each $\ell\in \La$, the associated effective volume is $v_\ell = \vor(\ell) $. For $\ell \in \L^\i$ the effective volume $v_\ell$ will depend on the geometry of the interface (see \cite{PRE-ac.2dcorners}), let $\omega_\ell:=\displaystyle{\frac{|v_\ell|}{|{\rm vor} (\ell)|}}$ denote the volume ratio of $v_\ell$ with respect to $\vor$. For each element $T \in \Th$ we define the effective volume of $T$ by 
\begin{displaymath}
\omega_T := |T \setminus (\bigcup_{\ell \in \L^{\a}} {\rm vor}(\ell))\setminus(\bigcup_{\ell \in \L^{\i}} v^\i_\ell)|.
\end{displaymath}
We note that $\omega_T =0 $ if $T\in \T^\a_h\setminus\T^\i_h$, $\omega_T= |T|$ if  $T\in \T^\c_h\setminus \T^\i_h$, and $0\leq \omega_T < |T|$ if $T\in \T^\i_h$.
The choices of $v_\ell$ and $\omega_T$ satisfy $\sum_{\ell\in \L^{\a,\i}} v_\ell + \sum_{T\in \Th} \omega_T = |\Omega_{h,R}|$.
	
Now we are ready to define the generic a/c coupling energy functional $\Eh$,
\begin{align}
  \label{eq:generic_ac_energy}
  \Eh(y_h) := & \sum_{ \ell \in \L^\a}  V_\ell(y_h)  + \sum_{\ell \in \L^\i} \omega_\ell V^\i_\ell(y_h)  + \sum_{T \in \Th} \omega_T  W(\D y_h|_T)   %\\
%  		 = & \sum_{\ell \in \L} V^h_\ell(\Ia u_h)  + \sum_{T \in \Th} \omega_T  W(\D u_h|_T)  .
\end{align}
where $V_\ell^\i$ is a modified interface site potential which satisfies consistency conditions \eqref{eq:force_pt} and \eqref{eq:energy_pt}. $\omega_\ell$ and $\omega_T$ are suitable coefficients, and their construction will be discussed immediately in Section \S~\ref{sec:formulation:grac} and references therein.

The goal of a/c coupling is to find  
\begin{equation}
  \label{eq:min_ac}
  y_{h,R} \in \arg\min \b\{ \Eh(y_h) \bsep y_h - y^B  \in
  \Us_{h,R} \b\}.
\end{equation}
The subscript $R$ in $y_{h,R}$ and $\Us_{h,R}$ can be omitted if there is no confusion.

\subsection{Consistent Atomistic/Continuum Formulation}
\label{sec:formulation:grac}
The construction of the interface potential in \eqref{eq:min_ac} is the key for the formulation of atomistic/continuum coupling methods. In order to demonstrate the a posteriori error estimate for the generic a/c coupling methods, we shall restrict ourselves to the GRAC type methods \cite{PRE-ac.2dcorners}. 

\subsubsection{The patch tests and consistent a/c method}
A key condition that has been widely discussed in the a/c coupling
literature is that $\Eh$ should exhibit no ``ghost forces''. We call this condition the \textit{force patch test}, namely, for $\L = \L^{\textrm{hom}}$ and $\Phi_{\ell} = \Phi$, 
\begin{equation}
  \label{eq:force_pt}
  \< \del \Eh(y^{\mF}), v \> = 0 \qquad \forall v \in \UsT,
  \quad \mF \in \R^{m \times d}.
\end{equation}
In addition, to guarantee that $\Eh$ approximates the atomistic
energy $\Ea$, it is reasonable to require that the interface
potentials satisfy an \textit{energy patch test}
\begin{equation}
  \label{eq:energy_pt}
  V_\ell^\i(y^{\mF}) = V(y^{\mF}) \qquad \quad \forall\ell \in \L^\i, \quad \mF \in \R^{m
    \times d}.
\end{equation}
If an a/c method satisfies the patch test \eqref{eq:force_pt} and \eqref{eq:energy_pt}, it is called a \textit{consistent a/c method}.

\subsubsection{GRAC: Geometric reconstruction based consistent a/c method}

To complete the construction of the consistent a/c coupling energy
\eqref{eq:generic_ac_energy}, we must specify the interface region $\L^\i$ and
the interface site potential. The geometric reconstruction approach was pioneered by Shimokawa
{\it et al} \cite{Shimokawa:2004}, and then modified and extended in
\cite{E:2006,PRE-ac.2dcorners}. We refer to \cite{COLZ2013} for details of the implementation of geometric reconstruction based consistent atomistic/continuum (GRAC) coupling energy for multibody potentials with general interaction range and arbitrary interfaces. The extension of GRAC to 3D is a work in progress \cite{Fang:2017}.

For a prototype implementation of GRAC, we consider the 2D triangular lattice $\Lhom :=\mA\Z^{2}$ with 
\begin{equation}
\mA = \mymat{1 & \cos(\pi/3) \\ 0 & \sin(\pi/3)}.
\label{eq:Atrilattice}
\end{equation}

Let $a_1 = (1,0)^T$, then $a_j = \mA_6^{j-1}a_1$, $j=1,\dots, 6$, are the nearest neighbour directions in $\Lhom$, where $\mA_6$ is the rotation matrix corresponding to a $\pi/3$ clockwise planar rotation. 

Given the homogeneous site potential $V\b( D y(\ell) \b)$, we can represent $V_\ell^\i$ in terms of $V$.  For each $\ell \in
\L^\i, \rho, \vsig \in \Rg_\ell$, let $C_{\ell;\rho,\vsig}$ be
free parameters, and define
\begin{equation}
  \label{eq:defn_Phi_int}
  V_\ell^\i(y) := V \B( \b( {\textstyle \sum_{\vsig \in
      \Rg_\ell} C_{\ell;\rho,\vsig} D_\vsig y(\ell) } \b)_{\rho \in
    \Rg_\ell} \B)
\end{equation}
A convenient short-hand notation is
\begin{displaymath}
  V_\ell^\i(y) = V( C_\ell \cdot Dy(\ell) ), \quad \text{where}
  \quad \cases{
    C_\ell := (C_{\ell;\rho,\vsig})_{\rho,\vsig \in \Rg_\ell}, \quad
    \text{and} &\\
    C_\ell \cdot Dy := \b( {\textstyle \sum_{\vsig \in
      \Rg_\ell} C_{\ell;\rho,\vsig} D_\vsig y } \b)_{\rho \in
    \Rg_\ell}. & }
\end{displaymath}

We name the parameters $C_{\ell;\rho,\vsig}$ as the \textit{reconstruction parameters}. They are chosen so that the resulting energy
functional $\Eh$ satisfies the energy and force patch tests
\eqref{eq:force_pt} and \eqref{eq:energy_pt}. A sufficient (and likely
necessary) condition for the energy patch test is that $\mF \cdot
\Rg_\ell = C_\ell \cdot (\mF \cdot \Rg)$ for all $\mF \in \R^{m
  \times d}$ and $\ell \in \L^\i$. This is equivalent to
\begin{equation}
  \label{eq:E_pt_grac}
  \rho = \sum_{\vsig \in \Rg_\ell} C_{\ell; \rho,\vsig} \vsig, \qquad
  \forall \ell \in \L^\i, \quad \rho \in \Rg_\ell.
\end{equation}
In addition, optimal condition and stabilisation mechanism were proposed in \cite{COLZ2013} and \cite{2013-stab.ac} to improve the accuracy and stability of GRAC scheme.

\def\Ta{\T_\a}
\def\Th{\T_{\rm h}}
\def\sh{\sigma^{\rm h}}

\subsubsection{Stress formulation}
\label{sec:formulation:stress}
The stress tensor based formulation can be obtained from the first variation of the energy. For any $y \in \Us$, and $y_h-y^B\in \Ush$, there exist piecewise constant tensor fields $\sa(y; \cdot)\in \PO(\Ta)^{2\times 2}, \sigc(y_h; \cdot) \in \PO(\Th)^{2\times 2}$, and $\sh(y_h; \cdot) \in \PO(\Th)^{2\times 2}$, such that they satisfy the following identities
\begin{align}
  \<\del\Ea(y),v\>&=\sum_{T\in\Ta}|T|\sa(y;T):\Dc{T} v, \forall v\in\Us, 
  \label{eq:atomstress}\\
  \<\del\Ec(y_h),v_h\>&=\sum_{T\in\Th}|T|\sigc(y_h;T):\Dc{T} v_h, \forall v_h\in\Ush,
  \label{eq:contstress}\\
    \<\del\Eh(y_h),v_h\>&=\sum_{T\in\Th}|T|\sh(y_h;T):\Dc{T} v_h, \forall v_h\in\Ush.
  \label{eq:acstress}
\end{align}
here $\Ta$ is the micro-triangulation induced by the reference lattice $\L$. 
We call $\sa$ an \textit{atomistic stress tensor}, $\sigc$ a \textit{continuum
stress tensor}, and $\sh$ an \textit{a/c stress tensor}. For the nearest neighbour interactions, we can choose the following atomistic stress tensor, continuum stress tensor, and a/c stress tensor respectively from the first variations \eqref{eq:atomstress}-\eqref{eq:acstress},
\begin{align}
    \label{eq:defn_Sa}
    \sa(y; T) :=~& \frac{1}{\det \mA} \sum_{b=(\ell, \ell+\rho)\in \partial T \bigcap \Bs} \partial_\rho V_\ell \otimes
    a_\rho, \\
    \label{eq:defn_Sc}
    \sigc(y_h; T) :=~& \pp W(\Dc{T} y_h)=\frac{1}{\det \mA} \sum_{j = 1}^6 \partial_j V(\nabla_T y_h) \otimes a_j, \\
    \label{eq:defn_Sh}
    \sh(y_h; T) :=~&\sum_{b=(\ell, \ell+\rho)\in \partial T \bigcap \Bs} \partial_\rho V_\ell^h(\Ia y_h) \otimes
    a_\rho + \omega_T \sigc(y_h; T).
\end{align}

We call piecewise constant tensor field $\sigma\in \PO(\T)^{2\times 2}$ \textit{divergence free} if 
\begin{displaymath}
	\sum_{T\in\T}|T|\sigma(T):\Dc{T} v\equiv 0, \forall v\in ({\rm P}_1(\T))^2.
\end{displaymath}
By definitions \eqref{eq:acstress}, it is easy to know that the force patch test condition \eqref{eq:force_pt} is equivalent to that $\sh(\mF x)$ is divergence free for any constant deformation gradient $\mF$.

\def\Fh {\mathcal{F}_h}

 The discrete divergence free tensor fields over the triangulation $\T$ can be characterized by the non-conforming Crouzeix-Raviart finite elements \cite{PRE-ac.2dcorners, Or:2011a}. The Crouzeix-Raviart finite element space over $\T$ is defined as
\begin{align*}
N_{1}(\T)=\{c:\bigcup_{T\in\T}\textrm{int}(T)\to\R \quad \big{|}& \quad c \textrm{ is piecewise affine w.r.t. }\T, \textrm{and} \\ 
 & \textrm{continuous in edge midpoints } q_{f}, \forall f\in\mathcal{F}\}
\end{align*}

The following lemma in \cite{PRE-ac.2dcorners} characterizes the discrete divergence-free tensor field.

\begin{lemma}
A tensor field $\sigma\in\mP_{0}(\T)^{2\times 2}$ is divergence free if and only if there exists a constant $\sigma_{0}\in\R^{2\times 2}$ and a function $c\in N_{1}(\T)^{2}$ such that 
\begin{displaymath}
% \begin{center}
\sigma = \sigma_{0} + \nabla c\mJ, \qquad \textrm{where}\quad\mJ = \mymat{0 & -1 \\ 1 & 0}\in{\sf SO}(2).
% \end{center}
\end{displaymath}
\label{lem:divfree}
\end{lemma}

The immediate corollary provides a representation of the stress tensor.

\begin{corollary}
\label{cor:divfree}
The stress tensors in the definitions \eqref{eq:atomstress}-\eqref{eq:acstress} are not unique. 
Given any stress tensor $\sigma \in \mP_{0}(\T)^{2\times 2}$ satisfies one of the definitions \eqref{eq:atomstress}-\eqref{eq:acstress} , where $\T$ is the corresponding triangulation. Define the admissible set as $\rm{Adm}(\sigma):=\{\sigma  + \nabla c \mJ, c\in N_1(\T)^2\}$, then any $\sigma'\in \rm{Adm}(\sigma)$  satisfies the definition of stress tensor.
\end{corollary}

\subsubsection{A Priori Error Estimates}

In the analytical framework proposed in \cite{LuOr:acta, 2013-defects}, the numerical error can be split into 3 parts: the \textit{modeling error} due to the discrepancy between the atomistic model and the continuum model at the interface and the finite element edges, the \textit{coarsening error} due to finite element discretization of the solution space in the continuum region, and the \textit{truncation error} due to the finite size of the computational domain. It is proven in \cite{2013-defects} that there exists a strongly stable solution $y_{h, R}$ to \eqref{eq:min_ac} and a constant $C^{\rm{a-priori}}$for GRAC method such that,

\begin{equation}
	\|\nabla u_{h,R} - \nabla u\|_{L^2} \leq C^{\rm{a-priori}} \big(\|hD^2 u\|_{\ell^2(\L\bigcap(\Omega^\c_R))}+\|Du\|_{\ell^2(\L\setminus B_{R/2})}\big)
\end{equation}
where $u_{h,R} = y_{h,R} - y^B$.
 
With the generic decay property \eqref{eq:ptdef:regularity}, and the following quasi-optimal conditions:
\begin{itemize}
\item the radius of the atomistic region $\T^\a_{h, R}$ satisfies,
\begin{equation}
	\underline{C} R_\a^{1+2/d}\leq R \leq \overline{C} R_\a^{1+2/d}, 
\end{equation}
\item
$\T^\c_{h, R}$ is a graded mesh so that the mesh size function $h(x) = \mathrm{diam}(T)$ for $x\in T\in \T^\c_{h, R}$ satisfies,
\begin{equation}
	|h(x)|\leq C^{\rm{mesh}}  \b(\frac{|x|}{R_\a}\b)^\beta, \text{  with  }1<\beta<\frac{d+2}{2}.
\end{equation}
\end{itemize}

It holds that there exists a constant $C_0>0$, depending on $C^{\rm{a-priori}}$, $\underline{C}$, $\overline{C}$, $C^{\rm{mesh}} $, and $\beta$ such that for $R$ sufficiently large,
\begin{equation}
	\|\nabla u_{h,R} - \nabla u\|_{L^2} \leq C_0 R^{-d/2-1}.
\end{equation}
In particular, when $d=2$, and when P1 finite elements are used in the continuum region, we have, 
\begin{equation}
	\|\nabla u_{h,R} - \nabla u\|_{L^2} \leq C_0 N^{-1},
\end{equation}
where $N$ is the overall degrees of freedom.

\section{Error Analysis}
\label{sec:error}

We present the a posteriori error analysis in this section. In \S~\ref{sec:error:residue}, we derive the residual estimate for the consistent GRAC a/c coupling scheme introduced in \S~\ref{sec:formulation:grac}. Then, we give a lower bound for the stability constant which is computable from the a/c solution $u_h$ in \S~\ref{sec:error:stability}. Finally, we put forward the main results Theorem \ref{thm:h1error} and Theorem \ref{thm:energyerror} in \S~\ref{sec:error:mainresults}.

\subsection{Residual Estimate}
\label{sec:error:residue}

To be more precise, we restrict ourselves to the case of nearest neighbour multibody interactions, namely, we use the so-called "grac23" method introduced in \cite{PRE-ac.2dcorners} as the a/c coupling mechanism. We will extend the formulation to general short-range multibody interactions in a future work and discuss it briefly in \S~\ref{sec:conclusion}.

\def\Pdef{\mathscr{P}^{\rm def}}
\def\Usc{\Us^c}
\def\Use{\Us^{1,2}}
\def\Usrh{\Us_R^h}

For lattice function $u: \L \to \R^m$, we denote its continuous and piecewise affine interpolant with respect to the micro-triangulation $\T_\a$
by $\Ia u$. Notice that $\L$ is a lattice with defect, we can construct the piecewise interpolant with respect to $\L^{\rm hom}$ by extending $u$ to vacancy sites, which will be introduced in \S~Appendix~\ref{sec:appendix:extension}. Identifying $u = \Ia u$, we can define the (piecewise constant) gradient $\nabla u = \nabla \Ia u: \R^m\to \R^{m\times d}$ and the spaces of compact and finite energy displacements, respectively, by 
%\begin{align*}
%	\Usc: & = \{u: \L\to \R^m | \supp(\nabla u)\text{ is compact}\}. \\
%	\Use : & = \{u:\L\to \R^m | \nabla u \in L^2\}.
%\end{align*} 
\begin{displaymath}
	\Usc := \{u: \L\to \R^m | \supp(\nabla u)\text{ is compact}\}.
\end{displaymath}
It can be shown that that $\Usc$ is dense in $\Use$ \cite{2013-defects}. 

The first variation of the atomistic variational problem \eqref{eq:min} is to find $y-y^B \in \Use$ such that 
\begin{equation}
\label{eqn:firstvariationea}
	\<\del\Ea(y), v\> = 0, \quad \forall v \in \Use.
\end{equation}

\def\Eh{\E^{\rm h}}
The first variation of the a/c coupling variational problem \eqref{eq:min_ac} is to find $y_h - y^B \in \Us_{h, R}$ such that
\begin{equation}
\label{eqn:firstvariationeh}
	\<\del\Eh(y_h), v_h\> = 0, \quad \forall v_h\in \Us_{h, R}.
\end{equation}

We introduce the truncation operator $T_R$ as in \cite{2013-defects} by first choosing a $C^1$ cut-off function $\eta(x) = 1$ for $|x|\leq 4/6$ and $\eta(x)=0$ for $|x|\geq 5/6$. Define $T_R: \Us^{1,2} \to \Us_R$  for $R>0$ by
\begin{displaymath}
	T_R u(\ell):= \eta(\ell/R)(u(\ell)-a_R), \text{  where  }a_R:=\int_{B_{5R/6}\setminus B_{4R/6}}I_a u(x)\dx,
\end{displaymath}
where $\Us_R$ is defined by
\begin{displaymath}
\Us_R:  = \{u \in \Usc | u(x) = 0 \ \forall x\in \L \backslash \Omega_R\}.
\end{displaymath}

The residual $\mR$ is defined as an operator on $\Use$ which is given by
\begin{equation}
\mR[v] = \< \del\Ea(\Ia y_h), v\>, \quad \forall v\in \Use.
\end{equation}
By \eqref{eqn:firstvariationeh}, denote $v_R = T_R v$, and take $v_h = \Cs_h T_R v: \Use\to \Us_{h,R} $, where $\Cs_h: \Us_R\to \Us_{h, R}$ is the modified Cl\'{e}ment operator \cite{clement:1975, Verf:1999} whose definition will be made clear in the following subsections. By \eqref{eqn:firstvariationeh} we can separate the residual into three groups,

\begin{align*}
	\mR[v] = \< \del\Ea(\Ia y_h), v\>  = & \< \del\Ea(\Ia y_h), v\> - \<\del\Eh(y_h), v_h\>\\
	 = & \< \del\Ea(\Ia y_h), v\> - \< \del\Ea(\Ia y_h), v_R\> \\
	& + \< \del\Ea(\Ia y_h), v_R\> - [ \del\Eh(y_h),  v_R] \\
	& + [ \del\Eh(y_h), v_R] - \<\del\Eh(y_h), v_h\>.
\end{align*}
Notice that $v_R\notin \Uhr$, therefore we cannot use the pairing $\<\del\Eh(y_h), v_R\>$. Instead, we define operation $[\cdot, \cdot]$ as,

\begin{align}
	[ \del\Eh(y_h), v_R]:= &\sum_{T\in\Th}\int_T\sh(y_h, T) \nabla v_R \dx \notag\\
				       = & \sum_{T\in\Th}\sh(y_h, T) (\sum_{T'\in\Ta, T'\bigcap T \neq \emptyset}|T \bigcap T'|\nabla v_R) \notag \\
				       = &\sum_{T\in \Ta}|T| \big(\sum_{T'\in \Th, T'\bigcap T\neq \emptyset}\frac{|T'\bigcap T|}{|T|}\sh(y_h, T')\big)\nabla v_R \label{eq:ehvr}
\end{align}

In the above decomposition of the residual $\mR[v]$, the first group $\mR_1:=\< \del\Ea(\Ia y_h), v\> -  \< \del\Ea(\Ia y_h), v_R\> $ 
represents the \textit{truncation error}, the second group $\mR_2:=\< \del\Ea(\Ia y_h), v_R\> - [ \del\Eh(y_h), v_R]$ represents the \textit{modeling error}, and the third group $\mR_3 := [ \del\Eh(y_h), v_R] -\<\del\Eh(y_h), v_h\>$ represents the \textit{coarsening error}.  We will deal with the contributions from those three groups separately in the following subsections.

\begin{remark}
\label{rem:nonunique}
Those residual estimators $\mR_1$, $\mR_2$ and $\mR_3$ are based on first variation of the energies, and can be in turn represented by stress formulation. By Lemma \ref{lem:divfree} and Corollary \ref{cor:divfree}, the stresses are unique up to a divergence-free tensor field. Therefore, we need to minimize those estimators with respect to divergence-free tensor field, which will be introduced in \S~\ref{sec:Interfacial-stress}.
\end{remark}

\subsubsection{Truncation error}
\label{sec:residue:truncation}

\def\Tr{\rm Tr}
To analyze the truncation error $\mR_1$, we need the Lemma 7.3 for the truncation operator $T_R$ in \cite{2013-defects}, namely, if the radius of the computational domain $R$ is sufficiently large (in the nearest neighbour case, we only need $R>6$), the following estimates hold
\begin{align*}
	\|Dv_R - Dv \|_{\ell^2} & \leq C^{\Tr} \|Dv\|_{\ell^2(\L\setminus B_{R/2})} \quad  \forall v \in \Us^{1,2}, \\
	\|Dv_R\|_{\ell^2} & \leq C^{\Tr} \|Dv\|_{\ell^2(\L\bigcap B_R)} \quad   \forall v \in \Us^{1,2},
\end{align*}
where $v_R = T_R v$, and $C^{\Tr}$ is independent of $R$. 
%\lz{change $R/2$ to $c_0 R$, as long as $(1-c_1)R\geq \rcut$}

For any $v\in\Use$, the stress-based formulation of the first variation \eqref{eq:stressform}, the fact that $ v_R(\ell) = v(\ell)$ for $|\ell/R|\leq 4/6$,  the equivalence of $\|D v \|_{\ell^2}$ and $\|\nabla v \|_{L^2}$, and Cauchy-Schwarz inequality lead to, 
\begin{align}
	|\mR_1| &= |\< \del\Ea(\Ia y_h), v\> - \< \del\Ea(\Ia y_h), v_R\>| \notag\\
	& = |\sum_{T\in \Ta} \sa(\Ia y_h, T) (\nabla v - \nabla v_R) - \sum_{T\in \Ta} \sigma^0 (\nabla v - \nabla v_R)| \label{eq:stressform}\\  
	& \leq \int_{\Omega_R \setminus B_{R/2}} |(\sa(\Ia y_h)-\sigma^0)(\nabla v - \nabla v_R)| \dx \notag\\
	& \leq \|\sa(\Ia y_h)-\sigma^0\|_{L^2(\Omega_R\setminus B_{R/2})}\|\nabla v - \nabla v_R\|_{L^2} \notag\\
	& \leq C^{\Tr}\|\sa(\Ia y_h)-\sigma^0\|_{L^2(\Omega_R \setminus B_{R/2})}\|\nabla v \|_{L^2} \label{eq:trunerr}
\end{align}
where $\sigma^0$ is divergence-free, i.e. $\sum_{T\in \Ta} \sigma^0 (\nabla v - \nabla v_R)=0$. In this paper, we assume a macroscopic applied strain $B\in\R^{d\times d}$, hence we can specify $\sigma^0 = \pp W(y^B)$. If we do not have uniform deformation at far field, for example in the case of nano-indentation, $\sigma^0$ can be computed from surface deformation.
Thus, the truncation error estimator $\eta_T$ is given by 
\begin{equation}
\label{eqn:etat}
	\eta_T(u_h): = C^{\Tr}\|\sa(\Ia u_h)-\sigma^0\|_{L^2(\Omega_R\setminus B_{R/2})}.
\end{equation}

\begin{remark}
\label{rem:etaT}
The numbers $4/6$, $5/6$ in the definition of truncation operator $T_R$, and consequently $R/2$ in the estimator $\eta_T$ are not essential. We can choose different numbers to define an estimator on a smaller outer domain, but the constant $C^{\Tr}$ will increase correspondingly. In practice, since $\T_h$ is a graded mesh, we can choose the boundary layer of triangles to evaluate $\eta_T$.
\end{remark}
%%%%%%%%%%%%%%%%%%%%%%%%%%%%%%%
\subsubsection{Modeling error}
\label{sec:residue:model}

In the analysis of the modeling error $\mR_2$, the stress based formulation of $\< \del\Ea(\Ia y_h), v_R\>$ and the definition of $[ \del\Eh(y_h), v_R]$ \eqref{eq:ehvr} lead to,
\begin{align}
	|\mR_2|:= & |\< \del\Ea(\Ia y_h), v_R\> - [ \del\Eh(y_h), v_R]| \nonumber\\
	        =&\big| \sum_{T\in \Ta}|T|\sa(\Ia y_h, T)\nabla v_R - \sum_{T\in \Ta}|T| \big(\sum_{T'\in \Th, T'\bigcap T\neq \emptyset}\frac{|T'\bigcap T|}{|T|}\sh(y_h, T')\big)\nabla v_R \big| \nonumber \\
		\leq& C^{\Tr}\big\{\sum_{T\in\Ta}|T|\big[\sa(\Ia y_h, T) - \sum_{T'\in \Th, T'\bigcap T\neq \emptyset}\frac{|T'\bigcap T|}{|T|}\sh(y_h, T')\big]^2\big\}^{\frac12}\|\nabla v\|_{L^2}.
	\label{eqn:error}
\end{align}
As a result, we define the modeling error estimator $\eta_M$ by,
\begin{equation}
\label{eqn:etam}
	\eta_M(y_h) := C^{\Tr}\big\{\sum_{T\in\Ta}|T|\big[\sa(\Ia y_h, T) - \sum_{T'\in \Th, T'\bigcap T\neq \emptyset}\frac{|T'\bigcap T|}{|T|}\sh(y_h, T')\big]^2\big\}^{\frac12}.
\end{equation}

With the canonical choice of $\sa$ and $\sc$ in \eqref{eq:defn_Sa} and \eqref{eq:defn_Sc}, we can see that only those $T\in \Ta$ intersects with the interface and edges in $\T^\c_h$ have nontrivial contributions to $\eta_M$. 
%%%%%%%%%%%%%%%%%%%%%%%%%%%%%%%
\subsubsection{Coarsening  error}
\label{sec:residue:coarse-graining}
For the coarsening error $\mR_3$, we first observe that
\begin{align}
\mR_3 := &[ \del\Eh(y_h), v_R] -\<\del\Eh(y_h), v_h\>, \nonumber \\
	=& \sum_{T\in\Th}  \int_T \sh(y_h, T) (\nabla v_R - \nabla v_h) \dx. \label{eqn:R3}
\end{align}

Here, we take $v_h = \Cs_h v_R$, where $\Cs_h$ is the modified Cl\'{e}ment interpolation operator \cite{clement:1975, Verf:1999}. For any node $x\in \Nh$ in the triangulation $\Th$, let $\phi_x$ be the nodal basis with respect to $x$ on $\Th$, and $\omega_x = \supp (\phi_x)$ be the support of $\phi_x$.  The interpolation operator $\Cs_h: L^1(\Omega_{h,R})\to V_h$ can be defined by, 
\begin{displaymath}
	\Cs_h w = \sum_{x\in\Nh \bigcap {\rm Int}(\Omega_h)} w_x\phi_x, \quad\text{where }\ w_x = \frac{\int_{\omega_x}w\phi_x\dx}{\int_{\omega_x} \phi_x \dx}, \forall x\in \Nh.
\end{displaymath}

By definition, $\Cs_h w$ satisfies the Dirichlet boundary condition. The Clement interpolation enjoys the following properties \cite{CarVer:1999, Verf:1999}, for any element $T\in \Th$, and any interior edge $f\in \Fh\bigcap {\rm int(\Omega_{h,R})}$,  
\begin{align}
	\|w- \Cs_h w\|_{L^2(T)} &\leq C_{\Th}h_T\|\nabla w\|_{L^2(\omega(T))},\label{eqn:clement1}\\
	\|w- \Cs_h w\|_{L^2(f)} & \leq C'_{\Th}h_f^\frac12\|\nabla w\|_{L^2(\omega(f))},\label{eqn:clement2}	
\end{align}
where $h_T$ is the diameter of $T$, and $h_f$ is the length of $f$. The element patch is $\omega(T) := \bigcup_{x\in \Nh\bigcap T}\omega_x$, and the edge patch is $\omega(f) := \bigcup_{x\in \Nh\bigcap f}\omega_x$. The constants $C_{\Th}$  and $C'_{\Th}$ depend 
only on the shape regularity of $\Th$. 

%We also have the stability of the Cl\'{e}ment interpolation, 
%\begin{displaymath}
%	\|\nabla \Cs_h w\|_{L^2(\Omega_R)} \leq C''_{\Th}\|\nabla w\|_{L^2(\Omega_R)}.
%\end{displaymath}
%where the constant $C''_{\Th}$ only depends on the shape regularity of $\Th$. 

\def\sjump{\llbracket\sh\rrbracket_f}
\def\dx {\rm dx}
\def\Fh {\mathcal{F}_h}
\def\Fhi{\Fh\bigcap {\rm int}(\Omega_R)}
For notational convenience, we assume that each interior edge $f\in \Fh\bigcap {\rm int}(\Omega_h)$ has a prescribed orientation. $\Tp_f$ and $\Tm_f$ are the triangles on the left hand side and right hand side of the edge $f$, $\np$ and $\nm$ are the corresponding outward unit norm vector. The integration by parts of \eqref{eqn:R3} leads to,
\begin{align*}
\mR_3  = &\sum_{T\in\Th}\int_T \sh(y_h, T) (\nabla v_R - \nabla v_h)\dx\\
	= & \sum_{f\in\Fh\bigcap {\rm int}(\Omega_R)} \int_f (\sh(y_h, \Tp_f)\np+\sh(y_h, \Tm_f)\nm)\cdot(v_R-v_h)\ds\\
	= & \sum_{f\in\Fh\bigcap {\rm int}(\Omega_R)} \llbracket\sh\rrbracket_f \cdot\int_{f\in\Fh} (v_R-v_h)\ds,\\
\end{align*}
where $\llbracket\sh\rrbracket_f:=\sh(y_h, \Tp_f)\np+\sh(y_h, \Tm_f)\nm$ denotes the jump of $\sh$ across the edge $f$. Cauchy-Schwarz inequality and the property of Clement interpolation \eqref{eqn:clement2} give rise to,
\begin{align*}
|\mR_3| & \leq \sum_{f\in\Fhi} |\sjump| h_f^\frac12 \|v_R-v_h\|_{L^2(f)}\\
			& \leq C'_{\Th}\sum_{f\in\Fhi} |\sjump| h_f \|\nabla v_R-\nabla v_h\|_{L^2(\omega_f)}\\
			& \leq C'_{\Th}(\sum_{f\in\Fhi} (h_f\sjump)^2)^\frac12  (\sum_{f\in\Fhi} \|\nabla v_R-\nabla v_h\|^2_{L^2(\omega_f)})^\frac12\\
			& \leq \sqrt{3}C'_{\Th}(\sum_{f\in\Fhi} (h_f\sjump)^2)^\frac12 \|\nabla v_R-\nabla v_h\|_{L^2(\Omega)}\\
  		& \leq \sqrt{3}C^{\Tr}C'_{\Th}(\sum_{f\in\Fhi} (h_f\sjump)^2)^\frac12 \|\nabla v\|_{L^2(\Omega)}.
\end{align*}
The coarse-graining error estimator is then defined as,
\begin{equation}
\label{eqn:etac}
\eta_C(u_h) := \sqrt{3}C^{\Tr}C'_{\Th}(\sum_{f\in\Fh} (h_f\sjump)^2)^\frac12
\end{equation}

\subsubsection{Residual Estimate}
\label{sec:residue:total}
Combining the above estimates, we have the following theorem for the residual.
\begin{theorem}
\label{thm:residue}
For $\forall v\in \Use$, let $y_h$ be the a/c solution of variational problem \eqref{eq:min_ac}, the residual  $\mR[v] = \< \del\Ea(\Ia y_h), v \>$ can be bounded by the sum of the truncation error (the $L^2$ norm of the atomistc stress tensor close to the outer boundary), modeling error (the difference of a/c stress tensor and atomistic stress tensor), and the coarsening error (jump of a/c stress tensor across interior edges), namely,
	\begin{equation}
		 \< \del\Ea(\Ia y_h), v \> \leq  \big(\eta_T(y_h) + \eta_M(y_h) + \eta_C(y_h) \big)\|\nabla v\|_{L^2},
		\label{eq:residue}
	\end{equation}
where $\eta_T(y_h) $, $\eta_M(y_h)$ and $ \eta_C(y_h)$ are given in \eqref{eqn:etat}, \eqref{eqn:etam} and \eqref{eqn:etac} respectively.	
\end{theorem}

\begin{remark}
\label{rem:stressform}
All the estimators $\eta_T$, $\eta_M$ and $\eta_C$ depend on the a/c solution $y_h$, through their dependence on the discrete stress tensor $\sh(y_h)$ and $\sa(\Ia y_h)$. We can therefore write,  
\begin{equation}
	\eta(y_h) := \tilde{\eta}(\sa(\Ia y_h), \sh(y_h)) = \eta_T(y_h) + \eta_M(y_h) + \eta_C(y_h).
	\label{eq:residue total}
\end{equation}

By Remark \ref{rem:nonunique} we denote ${\rm Adm}(\sh), {\rm Adm}(\sa)$ the sets of all possible stress tensors. Therefore, the desired estimate of the residual is
\begin{equation}
		 \< \del\Ea(\Ia y_h), v \> \leq  \min_{{\rm Adm}(\sh(y_h)), {\rm Adm}(\sa(\Ia y_h))}\tilde{\eta}(\sa(\Ia y_h), \sh(y_h))\|\nabla v\|_{L^2}.
		\label{eq:sharpresidue}
\end{equation}

We refer to the exact or approximate minimization of the residual with respect to the admissible tensor field as ``stress tensor correction", and we will discuss the implementation of stress tensor correction in detail in \S~\ref{sec:Interfacial-stress}.
\end{remark}

%%%%%%%%%%%%%%%%%%%%%%%%%%%%%%%%%%%%%%%%%%%%%%%%%%%%
%STABILITY%
%%%%%%%%%%%%%%%%%%%%%%%%%%%%%%%%%%%%%%%%%%%%%%%%%%%%
\subsection{Stability}
\label{sec:error:stability}

In this subsection, we will deduce a computable estimate of the a posteriori stability constant. Similar as the residual estimate, we restrict ourselves to the case of nearest-neighbour interaction with vacancies. We follow the stability analysis in \cite{OrtnerShapeev:2011}. The main difference is: first, we derive the stability results for the many-body potentials of generic pair functional form \eqref{eq:pairfunctional}, while in  \cite{OrtnerShapeev:2011} only pair interaction potentials are considered; second, in the a posteriori analysis the stability constant depends on the atomistic Hessian $\ddel\Ea$ and the a/c solution $u_h$, and therefore it is computable, as opposed to the a priori analysis in \cite{OrtnerShapeev:2011}, the stability constant is related to the a/c Hessian $\ddel\Eh$ and the unknown atomistic solution $u$ where certain assumptions for $u$ have to be made.

\begin{theorem}
\label{thm:stability}
Suppose that the multi-body interaction potential is of the generic pair functional form \eqref{eq:pairfunctional}, we have the following results,
\begin{equation}
\label{eq:Stability_A}
\<\delta^2\Ea(I_a y_{h})v,v\> \ge \gamma(y_h)  \| \nabla v \|_{L^2(\Omega_R)}^2 \quad
\forall v \in \Us ,
\end{equation}
where the precise definition of $\gamma(y_h)$ will be given as the analysis proceeds.
\end{theorem}

The proof of Theorem \ref{thm:stability} can be divided into the following steps: 

\begin{enumerate}[label=\theenumi.]
\item Write $\ddel\Ea(\Ia y_h)$ as a quadratic form with nonuniform coefficients defined on the interaction bonds;
\item Use the perturbation arguments \eqref{eq:Perturb_ineq_dot}, \eqref{eq:Perturb_ineq_cross} to bound $\ddel\Ea$ by quantities from a uniform deformation; 
\item Define the so-called vacancy stability index \eqref{eqn:vacancyindex} to further bound $\ddel\Ea$ for lattice with defects by the stability constant for a uniformly deformed homogeneous lattice; 
\item The stability constant can be obtained through an optimization procedure.
\end{enumerate}

Recall that by \eqref{eq:bonds}, $\Bs$ is the collection of all the nearest neighbour bonds in the reference lattice $\L$. Here we define  
\begin{equation}
\label{def: definition of bonds including vacancies}
\bbB:= \{ (\ell,\ell+\rho): \ell \in \Lhom, \rho \in \Rg_\ell\} 
\end{equation}
to be the collection of all the nearest neighbour bonds in the homogeneous reference lattice $\Lhom$. To simplify notation, we use $y$ to denote $I_a y_h$, and $\Omega$ to denote $\Omega_R$ in the following analysis of this section.

\subsubsection{Second variation of the energy}
Using the generic pair functional form multi-body interaction potential \eqref{eq:pairfunctional} and Remark \ref{rem:pairfuncational},  we write out the second variation of the atomistic energy $\Ea(y) = \sum_{\ell\in\L}V(|Dy(\ell)|)$ as
\begin{align}
\<\delta^2\Ea(y)v,v\> 
= &\sum_{\ell \in \L} \sum_{\rho,\vsig \in \Rg_{\ell}}
\partial_{\rho\vsig} V(|Dy(\ell)|) (D_\rho v(\ell))^{T} \b(\frac{D_\rho y(\ell)}{|D_{\rho} y(\ell)|} \otimes
\frac{D_{\vsig} y(\ell)}{|D_{\vsig} y(\ell)|} \b) (D_{\vsig} v(\ell)) \nonumber \\
%%%%%%%%%%%%%%%%%%%%%%%%%%%%%%%%%%%%%%%%%%
&+ \sum_{\ell \in \L} \sum_{\rho \in \Rg_{\ell}} \frac{\partial_{\rho} V(|Dy(\ell)|)}{|D_{\rho} y(\ell)|}
(D_{\rho} v(\ell))^{T} \b(\mI - \frac{D_{\rho} y(\ell)}{|D_{\rho} y(\ell)|} \otimes \frac{D_{\rho}
  y(\ell)}{|D_{\rho} y(\ell)|} \b) D_{\rho} v(\ell) \nonumber \\ 
%%%%%%%%%%%%%%%%%%%%%%%%%%%%%%%%%%%%%%%%%%
=& \sum_{\ell \in \L} \sum_{\rho \in \Rg_{\ell}}
\frac{\partial_{\rho\rho} V(|Dy(\ell)|)}{|D_{\rho} y(\ell)|^2} (D_{\rho} y(\ell) \cdot D_{\rho} v(\ell))^2 
\nonumber\\ 
&+\sum_{\ell \in \L } \sum_{\rho,\vsig \in \Rg_{\ell}, \rho \neq \vsig} \frac{\partial_{\rho\vsig}
  V(|Dy(\ell)|)}{|D_{\rho} y(\ell)| |D_{\vsig} y(\ell)|} (D_{\rho} y(\ell)\cdot D_{\rho} v(\ell))(D_{\vsig} y(\ell) \cdot D_{\vsig} v(\ell)) \nonumber\\ 
&+ \sum_{\ell \in \L} \sum_{\rho \in \Rg_{\ell}}
\frac{\partial_{\rho} V(|Dy(\ell)|)}{|D_{\rho} y(\ell)|^3} | D_{\rho} y(\ell) \times D_{\rho} v(\ell)|^2,
%%%%%%%%%%%%%%%%%%%%%%%%%%%%%%%%%%%%%%%%%%
\label{eq:Second_Vari_A}
\end{align}
where $\partial_{\rho} V(|Dy(\ell)|)$ represents the first order partial derivatives of $V(|Dy(\ell)|)$ with respect to $|D_{\rho}y(x)|$, and $\partial_{\rho\vsig} V(|Dy(\ell)|)$ represents the second order partial derivatives with respect to $|D_{\rho}y(\ell)|$ and $|D_{\vsig} y(\ell)|$, 
$\mI$ is the identity matrix, and $a \times b =
a_1b_2-a_2b_1$. We have also used the identity
\begin{align}
h_1^T(\frac{r_1}{|r_1|} \otimes \frac{r_2}{|r_2|})h_2 &= (h_1 \cdot
\frac{r_1}{|r_1|}) (h_2 \cdot \frac{r_2}{|r_2|}),  \nonumber \\
\text{ and } h^T(\mI-\frac{r}{|r|} \otimes \frac{r}{|r|})h &= |h \times
\frac{r}{|r|}|^2.
\label{eq:Second_Vari_id_2}
\end{align}

For nearest neighbour interactions, $|\Rg(\ell)| \le 6$, we define 
\begin{align*}
C^1_{\ell,\rho} &=   \frac{\partial_{\rho\rho}V(D y(\ell))}{|D_{\rho}  y(\ell)|^2},\quad
 C^2_{\ell,\rho} =  0\wedge \min_{\vsig,\vsig \neq \rho} \frac{\partial_{{\rho}{\vsig}}
   V(D y(\ell))}{|D_{\rho} y(\ell)||D_{\vsig} y(\ell)|} ,\\
C_{\ell,\rho} &= \min_{\ell}(C^1_{\ell,\rho}  - 5C^2_{\ell,\rho}), \quad 
C^{\perp}_{\ell,\rho} =  \frac{\partial_\rho V(Dy(\ell))}{|D_{\rho} y(\ell)|^3}.
\end{align*}
Applying Cauchy-Schwarz inequality to \eqref{eq:Second_Vari_A}, we obtain the following estimate, 
\begin{align}
\<\delta^2\Ea(y)v,v\> \ge &  \sum_{\ell \in \L} \sum_{\rho\in\Rg_\ell} C_{\ell, \rho} |D_{\rho} y(\ell) \cdot D_{\rho} v(\ell)|^2 + \sum_{\ell \in \L} \sum_{\rho \in \Rg_\ell } C^{\perp}_{\ell, \rho} |D_{\rho} y(\ell) \times D_{\rho} v|^2 \nonumber \\
=& \sum_{b \in \Bs} C_b|D_b y(\ell) \cdot D_b v(\ell)|^2 + 
 \sum_{b \in \Bs} C^{\perp}_b|D_b y(\ell) \times D_b v(\ell)|^2\nonumber\\
 \geq & C\sum_{b \in \Bs} |D_b y(\ell) \cdot D_b v(\ell)|^2 + 
 C^{\perp}\sum_{b \in \Bs} |D_b y(\ell) \times D_b v(\ell)|^2\nonumber\\
=& C\sum_{b \in \Bs} \int_b |D_b y\cdot \nabla_b v|^2 \db + 
C^{\perp}\sum_{b \in \Bs}  \int_b |D_b y \times \nabla_b v|^2 \db.
\label{eq:Stab_1st_ineq}
\end{align}
where $C_b := C_{\ell, \rho}$ and $C_b^{\perp}:=C_{\ell, \rho}^{\perp}$ for $b=(\ell, \ell+\rho)$, $C^{(\perp)}: = \min_{b\in\Bs}C_b^{(\perp)}$ (here we use $C^{(\perp)}$ to denote both $C$ and $C^{\perp}$ for brevity). We have also used the fact that for nearest neighbour interactions, $D_b v = \nabla_b v(x)$, $\forall x \in \rm{int}(b)$, and $D_b y = D_b y(\ell)$ is a constant for each $b = (\ell, \ell+\rho) \in \Bs$. 

\subsubsection{The perturbation argument}

Our next task is to obtain the estimates, 
\begin{equation}
\label{eq: stability step 2 overview}
C\sum_{b \in \Bs} \int_b |D_b y \cdot \nabla_b v|^2 \db \ge c
\|\nabla v\|^2_{L^2(\Omega)},
\text{ and } 
C^{\perp}\sum_{b \in \Bs} \int_b |D_b y \times \nabla_b v|^2 \db \ge c^{\perp}
\|\nabla v\|^2_{L^2(\Omega)},
\end{equation}
for some $c>0$ and  $c^{\perp}$ (which could be negative). 

\eqref{eq: stability step 2 overview} is not straighforward
since $D_b y$ varies on each $b\in \Bs$. To tackle this issue, we use the following perturbation results from Lemma 6.3 of \cite{OrtnerShapeev:2011}.
For $g \in \R^2$, $b \in \Bs$, and $\alpha > 0$, we have

\begin{align}
\B| |D_b y \cdot g|^2 - |\mB {\rho}_b \cdot g|^2 \B| &\le \alpha |\mB {\rho}_b
\cdot g|^2 + (1+\frac{1}{\alpha}) \Delta^2|{\rho}_b |^2|\mB^{T} g|^2,
\label{eq:Perturb_ineq_dot}  \\
\text{ and } \B| |D_b y \times g|^2 - |\mB {\rho}_b \times g|^2 \B| & \le \alpha^\perp |\mB {\rho}_b
\times g|^2 + (1+\frac{1}{\alpha^\perp}) \Delta^2|{\rho}_b|^2|\mB^{T} g^{\perp}|^2.
\label{eq:Perturb_ineq_cross} 
\end{align}
where $\rho_b$ is the direction vector of $b$, $\mB \in \R^{2 \times 2}$ is fixed, $\alpha^{(\perp)}$ are unknowns to be determined, and $\Delta  = \max_{T\in\T} \|\mB^{-1} \nabla y|_{T} - \mI\|$, $g^{\perp}$ is obtained by $\pi/2$ counterclockwise rotation of $g$. 

\def\bB{\mathbb{B}}

Given $y$, $\Delta$ and $\mB$ can be solved from the convex optimization problem $\Delta  = \max_{T\in\T} \|\mB^{-1} \nabla y|_{T} - \mI\|$. We will choose free parameters $\alpha$ and $\alpha^{\perp}$ in the subsequent analysis to keep the estimate of the stability constant sharp. Applying\eqref{eq:Perturb_ineq_dot} and \eqref{eq:Perturb_ineq_cross} to \eqref{eq:Stab_1st_ineq}, taking the same $\alpha$ and $\alpha^{\perp}$ for each bond $b\in \Bs$ and using the fact that $|{\rho}_b| = 1$ , we obtain
\begin{align}
&C \sum_{b \in \Bs} \int_b |D_b y\cdot \nabla_b v|^2 \db + C^{\perp}
\sum_{b \in \Bs}  \int_b |D_b y\times \nabla_b v|^2 \db \nonumber \\
%%%% line 2 %%%%%
\ge & C \sum_{b \in \Bs} \int_b |\mB {\rho}_b\cdot \nabla_b v|^2 \db 
+ C^{\perp} \sum_{b \in \Bs}  \int_b |\mB {\rho}_b\times \nabla_b v|^2 \db \nonumber \\
%%%%%%%%%%%%
&- \bigg( \alpha |C|\sum_{b \in \Bs} \int_b   |\mB {\rho}_b\cdot \nabla_b v|^2 \db
+ \alpha^{\perp} |C^{\perp}|\sum_{b \in \Bs} \int_b   |\mB {\rho}_b\times \nabla_b v|^2 \db\nonumber\\
%%%%%%%%%%%%
& +\Delta^2 C(1+\frac{1}{\alpha})\sum_{b\in\bB} \int_b   |\mB^{T} \nabla_b v|^2 \db + 
\Delta^2C^{\perp} (1+\frac{1}{\alpha^{\perp}}) \sum_{b\in\bB}  \int_b   |\mB^{T} \nabla_b v^{\perp}|^2 \db \bigg) \nonumber \\
%%%%%%%%%%%%
= &   \tilde{C}\sum_{b \in \Bs} \int_b  |{\rho}_b\cdot \nabla_b v_{\mB}|^2 \db + 
 \tilde{C}^{\perp}\sum_{b \in \Bs} \int_b  |{\rho}_b\times \nabla_b v_{\mB}|^2 \db \nonumber\\
& -\bigg(\Delta^2C(1+\frac{1}{\alpha})\sum_{b\in\bB} \int_b |\mB^{T} \nabla_b Ev|^2 \db + 
\Delta^2C^{\perp} (1+\frac{1}{\alpha^{\perp}})\sum_{b\in\bB} \int_b   |\mB^{T} \nabla_b Ev^{\perp}|^2 \db \bigg) \nonumber 
\label{eq:Stab_2nd_ineq}
\end{align}
where $\tilde{C}^{(\perp)}:=C^{(\perp)} - \alpha |C^{(\perp)}|$, we have used $\mB^{T} \nabla_b v = \nabla_b \mB^{T} v$, $\mB {\rho}_b \cdot \nabla_b v = {\rho}_b \cdot \mB^T \nabla_b v$, and $v_{\mB} := \mB^T v$. $Ev$ is the extension of $v$ from $\L$ to the vacancy sites defined in the Appendix \S~\ref{sec:appendix:extension}, it is clear that $Ev^\perp = (Ev)^\perp$.

\def\cL{\mathcal{L}}

Let 
\begin{equation}
	\<\tilde{H}v, v\> :=  \tilde{C}\sum_{b \in \Bs} \int_b  |{\rho}_b\cdot \nabla_b v_{\mB}|^2 \db + 
\tilde{C}^{\perp} \sum_{b \in \Bs} \int_b  |{\rho}_b\times \nabla_b v_{\mB}|^2 \db \nonumber 
\end{equation}
and 
\begin{align}
	\<\tilde{\cL}^{(\perp)}v, v\> :&= 
C^{(\perp)} (1+\frac{1}{\alpha^{(\perp)}})\sum_{b\in\bB}   \int_b   |\mB^{T} \nabla_b Ev^{(\perp)}|^2 \db\nonumber\\
						 &= \tilde{L}^{(\perp)}\|\nabla (\mB^TEv^{(\perp)})\|^2_{L^2(\Omega)}.
						 \label{eq:tildeL}	
\end{align}
where $\displaystyle\tilde{L}^{(\perp)} = \frac{3}{\det \mA_6}(1+\frac{1}{\alpha^{(\perp)}})C^{(\perp)}$. 
\eqref{eq:tildeL} is due to the application of the so-called bond-density lemma with respect to Dirichlet boundary conditions \cite[Lemma 4.5]{Shapeev:2010a}. Combining the above results, we have the following estimate,
\begin{equation}
\label{eq:hessian}
	\<\ddel \Ea(y) v, v\>\geq \<\tilde{H}(y)v, v\> - \Delta^2(\tilde{L}\|\nabla (\mB^TEv)\|^2_{L^2(\Omega)}+\tilde{L}^{\perp}\|\nabla (\mB^TEv^{\perp})\|^2_{L^2(\Omega)})
\end{equation}

\subsubsection{Vacancy stability index}
We introduce the vacancy stability index $\kappa$ as
\begin{equation}
\kappa(\bbV) = \max\B\{ k >0: \Phi_{\Bs} (u) \ge k \Phi_{\bbB} (Eu), \ \forall
u \in \Us\B\}.
\label{eqn:vacancyindex}
\end{equation}
Since $\tilde{C}>0$ and $\tilde{C}^{\perp}$ might be negative, we define the constants 
\begin{equation}
	\bar{C}^{(\perp)}:=\min(\tilde{C}^{(\perp)}, \kappa\tilde{C}^{(\perp)}).
\end{equation}
We can further estimate \eqref{eq:hessian} by 
\begin{align}
\<\ddel\Ea(y) v, v \> \ge &\bar{C}   \sum_{b \in \bbB} \int_b |{\rho}_b\cdot \nabla_b (E\mB^T v) |^2 \db
+\bar{C}^{\perp}  \sum_{b \in \bbB} \int_b |{\rho}_b\times \nabla_b (E\mB^T v)|^2 \db \nonumber\\
&- \Delta^2(\tilde{L}\|\nabla (\mB^TEv)\|^2_{L^2(\Omega)}+\tilde{L}^{\perp}\|\nabla (\mB^TEv^{\perp})\|^2_{L^2(\Omega)}).
\label{eq:atomhessian}
\end{align}

\subsubsection{Stability of the homogenous lattice}
Now we need the stability estimates for the homogeneous lattice. Let
\begin{equation}
\<\bar{\Hs} v, v \> = \bar{C}   \sum_{b \in \bbB} \int_b |{\rho}_b \cdot \nabla_b (E\mB^T v)|^2 \db+\bar{C}^{\perp}  \sum_{b \in \bbB} \int_b |{\rho}_b\times \nabla_b (E\mB^T v)|^2 \db.
\end{equation}

By Lemma 6.4 of \cite{OrtnerShapeev:2011}, we have
\begin{equation}
	\<\bar{\Hs}v, v\> \geq \bar{\gamma} \| \nabla E\mB^T v\|^2_{L^2(\Omega)}.
\end{equation}
where $\displaystyle \bar{\gamma}:=\min(\frac{3}{4}\bar{c}+\frac{9}{4}\bar{c}^{\perp}, \frac{9}{4}\bar{c}+\frac{3}{4}\bar{c}^{\perp})$, and $\displaystyle\bar{c}^{(\perp)} = \frac{3}{\det \mA} \bar{C}^{(\perp)}$.

Furthermore, by the inequality \eqref{eq:extensionineq} for the extension operator $E$ in the appendix, we can estimate the stability of atomistic Hessian \eqref{eq:atomhessian} by,
\begin{equation}
	\<\ddel \Ea(y)v, v\> \geq \gamma(y) \|\nabla v\|^2_{L^2(\Omega)}.
\end{equation}
where 
\begin{equation}
	\gamma(y) = \frac{1}{3}\|\mB^{-T}\|_F^{-1}\bar{\gamma} - \Delta^2\|\mB\|_F^2(\tilde{L} + \tilde{L}^\perp). 
\label{eq: stability constant}
\end{equation}
    
%$\gamma = \min(\gamma_1+\gamma_1^{\perp}, \gamma_2 +\gamma_2^{\perp})$.
%such that 
%\begin{align}
%	&\gamma_1 :=3/4\bar{c} - \Delta^2\tilde{L}, \quad \gamma_1^\perp : = \frac{9}{4}\bar{c}^\perp - \Delta^2 \tilde{L}^\perp,\\
%	&\gamma_2 :=9/4\bar{c} - \Delta^2\tilde{L}, \quad \gamma_1^\perp : = \frac{3}{4}\bar{c}^\perp - \Delta^2 \tilde{L}^\perp,
%\end{align}

\subsubsection{Numerical Justification}
Tracing back the derivation of the stability constant $\gamma$, the only free parameters are $\alpha$, $\alpha^\perp$. Consequently, we can find the optimal $\gamma$ by maximization with respect to $\alpha $ and $\alpha^\perp$. 

We justify our a posteriori estimate for the stability constant of the atomistic Hessian numerically. We apply the same EAM potential as in \S~\ref{sec:numerics:problem} and take isotropic stretch $\mathrm{S}$ and shear loading $\gamma_{II}$ by setting
\begin{displaymath}
{\sf B}=\left(
	\begin{array}{cc}
		1+\mathrm{S} & \gamma_{II} \\
		0            & 1+\mathrm{S}
	\end{array}	 \right)
	\cdot{\sf {F_{0}}},
\end{displaymath}
where ${\sf F_{0}} \propto \mathrm{I}$ minimizing the corresponding Cauchy-Born energy density $W(F)$. The numerical results are listed in the following tables, where $\lambda$ stands for the smallest eigenvalue of atomistic Hessian, and $\gamma$ represents the optimal estimate of the stability constant.

% table 1. identity.
\begin{table}[h]
\begin{center}
\begin{tabular}{|c|c|c|c|}
\hline
number of vacancies & 0 & 1 & 2 \\
\hline
$\lambda$ & 17.436 & 14.107 & 12.905 \\
\hline
$\gamma$ & 5.284 & 2.206 & 1.451 \\
\hline
\end{tabular}
\caption{In this example, we test the stability for the reference configuration, namely, $\mathrm{S} = \gamma_{II} = 0$. The degrees of freedom of the atomistic model is about $3\times10^{4}$.}
\label{table:identity}
\end{center}
\end{table}

% table 2. stretch and shear.
\begin{table}[h]
\begin{center}
\begin{tabular}{|c|c|c|c|}
\hline
number of vacancies & 0 & 1 & 2 \\
\hline
$\lambda$ & 11.125 & 9.809 & 8.946 \\
\hline
$\gamma$ & 3.159 & 0.468 & -0.258 \\
\hline
\end{tabular}
\caption{In this example, we test the stability for the deformed configuration with $\mathrm{S} = \gamma_{II} = 0.03$. The degrees of freedom of the atomistic model is about $3\times10^{4}$.}
\label{table:deform}
\end{center}
\end{table}
From the numerical results, our estimates indeed give lower bound of the minimal eigenvalue of atomistic Hessian, however, the estimate may become negative when the deformation and number of vacancy sites increase. 
%To compute $\Delta$, $\mB$ and then the stability constant, we used CVX, a package for specifying and solving convex programs \cite{cvx}.

%%%%%%%%%%%%%%%%%%%%%%%%%%%%%%%%%%%%%%%%%%%%%%%%%%%%
%A Posteriori Estimates%
%%%%%%%%%%%%%%%%%%%%%%%%%%%%%%%%%%%%%%%%%%%%%%%%%%%%
\subsection{Main results}
\label{sec:error:mainresults}

We present the main theorems for the a posteriori errors in $H^1$ norm and energy in this section. 

\subsubsection{A Posteriori Error Estimates in $H^1$ norm}
\label{sec:error:mainresults:h1error}

We will need the following quantitative version of the inverse function theorem in \cite{LuOr:acta}.

\begin{lemma}
\label{lem:ift}
Let $X$ be a Hilbert space, $w_0\in X$, $R$, $M > 0$, and $E\in C^2(B^X_R(\omega_0))$ with Lipschitz continuous Hessian, 
$\|\delta^2E(x)-\delta^2E(y)\|_{L(X, X^*)}\leq M \|x-y\|_X$ for $x, y\in B^X_R(\omega_0)$. Suppose, moreover, that there exists constants $c$, $r>0$, such that
\begin{equation}
\<\delta^2 E(w_0)v, v\>\geq c \|v\|^2_X,\quad, \|\delta E(w_0)\|_{X^*} \leq r,\quad \text{ and } 2Mrc^{-2} < 1.
\end{equation}
Then there exists a unique $\bar{w} \in B^X_{2rc^{-1}}(w_0)$ with $\delta E(\bar{w}) = 0$ and
\begin{displaymath}
\< \delta^2 E ( \bar{w} ) v , v \> \geq  (1 - 2 M r c^{-2} )  c \| v \|^2_X .
\end{displaymath}
\end{lemma}

Take $X = \Ush$, $\omega_0$ as the a/c solution $y_h$ of \eqref{eqn:firstvariationeh}, and $M$ as the Lipschitz constant of $\ddel \Ea$. Combine the residual estimate in Theorem \ref{thm:residue}, stability estimate in Theorem \ref{thm:stability}, and Lemma \ref{lem:ift}, we have the following theorem for the a posteriori existence and error estimate.

\begin{theorem}
\label{thm:h1error}
Let $y_h$ be the a/c solution of \eqref{eqn:firstvariationeh}, $\eta(y_h)$ be the residual defined in \eqref{eq:residue total}, $ \gamma(y_h)$ be the stability constant defined in \eqref{eq: stability constant}, and $M$ be the Lipschitz constant of $\ddel \Ea$. Under the assumption that $\gamma(y_h)>0$ and $2M\eta(y_h)< \gamma(y_h)^2$, there exists a unique $y$ satisfying $y - y^B \in \Use$ which solves the atomistic variational problem \eqref{eqn:firstvariationea}, and satisfies the following error bound, 
\begin{equation}
	\|\nabla \Ia y_h - \nabla y\|_{L^2}\leq 2\frac{\eta(y_h)}{\gamma(y_h)},
\end{equation}
and the strong stability condition,
\begin{equation}
	\< \delta^2 E (y) v , v \> \geq  \big(1 - 2 \frac{M \eta(y_h)}{ \gamma(y_h)^2} \big)  \gamma(y_h) \| \nabla v \|^2_{L^2}, \quad\forall v\in\Use.
\end{equation}

\end{theorem}

\begin{remark}
Alternatively, the a posteriori error estimate can be deduced by the following argument in \cite{OrtnerWang:2014}, but we need to assume the existence of the atomistic solution $y$ and the closeness of $y$ to $\Ia y_h$ in $W^{1, \infty}$. By mean value theorem, there exists $\theta\in {\rm conv}\{y, \Ia y_h\}$ such that 
\begin{align}
	\<\ddel \Ea(\theta)v, v\> & =  \< \del\Ea(\Ia y_h), v\> - \<\del\Ea(y), v\> \nonumber\\
				  & =  \< \del\Ea(\Ia y_h), v\> \nonumber\\
				  & \leq  \eta(y_h)\|\nabla v\|_{L^2(\Omega)}.
				  \label{eq:residue2}
\end{align}
%\item $\displaystyle\|\nabla u_h-\nabla u\|_{L^\infty} \leq \ \frac{c_a(u_h)}{2C_{Lip}}$.

Combining the coercivity of $\Ea$ at $\Ia y_h$,
\begin{displaymath}
	\<\ddel \Ea(\Ia y_h)v, v\> \geq \gamma(y_h) \|\nabla v\|^2_{L^2},
\end{displaymath}
and the Lipschitz continuity (Fr\'{e}chet differentiability) of $\ddel \Ea$, we obtain that
\begin{align}
	\<\ddel \Ea(\theta)v, v\>  & \geq \<\ddel \Ea(\Ia y_h)v, v\>  - M\|y-\Ia y_h\|_{W^{1,\infty}} \|\nabla v\|_{L^2}^2\nonumber\\
	& \geq ( \gamma(y_h) - M\|y-\Ia y_h\|_{W^{1,\infty}} ) \|\nabla v\|^2_{L^2}
	\label{eq:stability2}
\end{align}
Let $v=y-\Ia y_h$ in \eqref{eq:stability2}, using \eqref{eq:residue2}, we have
\begin{equation}
	\|\nabla y - \nabla \Ia y_h\|_{L^2}\leq \frac{2\eta(y_h)}{\gamma(y_h)}
\end{equation}
if the closeness assumption $\displaystyle\|\nabla y_h-\nabla y\|_{L^\infty} \leq \ \frac{\gamma(y_h)}{2M}$ holds true. 

\end{remark}

\subsubsection{A Posteriori Error Estimate for the Energy}
\label{sec:error:mainresults:energyerror}

Total energy is an important physical quantity to be approximated in applications. In this section, we will derive an estimate for the energy difference $\Ea(y)-\Eh(y_h)$. The energy difference can be split into the sum of 
$\Ea(y)-\Ea(\Ia y_h)$ and $\Ea(\Ia y_h) - \Eh(y_h)$, thus,
\begin{equation}
|\Ea(y)-\Eh(y_h)| \leq  |\Ea(y)-\Ea(\Ia y_h)| + |\Ea(\Ia y_h) - \Eh(y_h)|
\label{eq: energy error separation}
\end{equation}

\def\ds{\rm ds}
For the first part, since $\Ea$ is twice differentiable along the segment $\{(1-s)y+s\Ia y_h | s\in(0,1)\}$, we obtain,
\begin{align}
	|\Ea(y)-\Ea(\Ia y_h)| & = |\int^1_0 \< \del \Ea((1-s)y+s\Ia y_h), y - y_h \>\ds|\nonumber\\
				       & = |\int^1_0 \< \del \Ea((1-s)y+s\Ia y_h) - \del \Ea(y), y - \Ia y_h \>\ds|\nonumber\\
				       & \leq M \| Dy - D\Ia y_h\|^2_{\ell^2 }\nonumber\\
				       & \leq M \|\nabla y - \nabla \Ia y_h\|^2_{L^2}.
\end{align}	
which can be further estimated by Theorem \ref{thm:h1error}, the constant $M$ is the Lipschitz constant of $\ddel\Ea$ which is independent of $y_h$. 

For the second part, let $\mu_E(y_h):=\Ea(\Ia y_h) - \Eh(y_h)$. We can rewrite $\Ea$ in the site based form,  
\begin{displaymath}
	\Ea(\Ia y_h) = \sum_{T\in\Ta} \frac{1}{6} \sum_{\ell\in T\bigcap \L} V_\ell(\Ia y_h).
\end{displaymath}
Moreover,  given $\Eh$ of the form  \eqref{eq:generic_ac_energy}, assuming for simplicity $\omega^\i_\ell = 1$, and $\T_h^\i$ is a few layers of atomistic micro-triangulation around the $\T_h^\a$,  which is actually the case for the implementation in \cite{PRE-ac.2dcorners}, we can rewrite $\Eh$ as follows,
\begin{align*}
	\Eh(y_h) = & \sum_{T\in\T_h^\a} \frac{1}{6} \sum_{\ell\in T\bigcap \La} V_\ell(\Ia y_h) + \sum_{T\in\T_h^\a} \frac{1}{6} \sum_{\ell\in T\bigcap \Li} V^\i_\ell(\Ia y_h)  +\\
						 & \sum_{T\in\T_h^\i\bigcap\T_h^\c} \Big\{ \sum_{\ell\in T\bigcap \Li}\frac{1}{6} \sum_{\ell\in T} V^\i_\ell(\Ia y_h) + (1-\frac{\#\{\ell\in T\bigcap \Li\}}{3})|T|W(\D\Ia y_h) \Big\}+\\
						 & \sum_{T\in\Th^{\c}\setminus\T_h^\i}\sum_{T'\in\Ta, T'\bigcap T\neq \emptyset}|T\bigcap T'|W(\D\Ia y_h).
\end{align*}

Hence $\mu_E$ can be expanded as,
\begin{align}
	\mu_E(y_h) = & \sum_{T\in\T_h^\a} \frac{1}{6} \sum_{\ell\in T\bigcap \Li} \big(V_\ell(\Ia y_h) - V^\i_\ell(\Ia y_h) \big)+ \nonumber\\
	& \sum_{T\in\T_h^\i\bigcap\T_h^\c} \Big\{\sum_{\ell\in T\bigcap \L^\i}\frac{1}{6} \sum_{\ell\in T} V_\ell( \Ia y_h) - \sum_{\ell\in T\bigcap \L}\frac{1}{6} \sum_{\ell\in T} V^\i_{\ell}(\Ia y_h )+\nonumber\\
						 & (1-\frac{\#\{\ell\in T\bigcap \Li\}}{3})|T|W(\D\Ia y_h)\Big\}+\nonumber\\
						 & \sum_{T\in\T_h^{\c}\setminus \T_h^\i}\sum_{T'\in\Ta, T'\bigcap T\neq \emptyset}\frac{|T\bigcap T'|}{|T'|}(\frac{1}{6} \sum_{\ell\in T'} V_\ell(\Ia y_h)-W(\nabla y_h)).
						  \label{eq: energy estimate II}
\end{align}

We note that the summand in the last term, which is summed over  $T\in\Th^{\c}$, is nonzero only if $\omega(T')\bigcap \partial T \neq \emptyset$, therefore can be rewritten as 
\begin{displaymath}
\sum_{T\in\Th^{\c}}\sum_{T'\in\Ta, \omega(T')\bigcap \partial T\neq \emptyset}\frac{|T\bigcap T'|}{2|T'|}(\frac{1}{3} \sum_{\ell\in T} V(D \Ia y_h (\ell))-V(\D\Ia y_h\rho)),
\end{displaymath}
noticing that $V_{\ell} = V$ when $T\bigcap \La =\emptyset$.

Hence we have the following theorem, 

\begin{theorem}
\label{thm:energyerror}
Given the same conditions in Theorem \ref{thm:h1error}, the difference of the energy can be bounded by the following inequality, 
\begin{displaymath}
	|\Ea(y)-\Eh(y_h)| \leq C^E(\eta(y_h))^2 +  |\mu_E(y_h)|.
\end{displaymath}
where $\displaystyle C^E=\frac{4M}{\gamma(y_h)^2}$, $\eta(y_h)$ and $\mu_E(y_h)$ are defined in \eqref{eq:residue total} and \eqref{eq: energy estimate II} respectively.
\end{theorem}

We denote the energy estimator by
\begin{equation}
	\eta_E(y_h) := C^E(\eta(y_h))^2 +  |\mu_E(y_h)|.
	\label{eq:etaE}
\end{equation}

\section{Adaptive Algorithms and Numerical Experiments}
\label{sec:numerics}
\def\Ta{\T_\a}
\def\Th{\T_{\rm h}}
\def\sh{\sigma^{\rm h}}
\def\sjump{\llbracket\s\rrbracket}

\newcommand{\fig}[1]{Figure \ref{#1}}
\newcommand{\tab}[1]{Table \ref{#1}}

In this section, we propose an adaptive mesh refinement algorithm based on the a posteriori error estimates in Theorem \ref{thm:h1error} and Theorem \ref{thm:energyerror}. Numerical experiments show that our algorithm achieves an optimal convergence rate in terms of accuracy vs. the degrees of freedom, which is the same as the a priori error estimates.

\subsection{Adaptive mesh refinement algorithm.}
\label{sec:numerics:algorithm}

Our goal is to design adaptive refinement algorithms by utilizing the residual based error estimators $\eta_M$, $\eta_C$, $\eta_T$ in \S~\ref{sec:error:residue} and $\mu_E$ in \S~\ref{sec:error:mainresults:energyerror}. The algorithm follows the usual Solve-Estimate-Mark-Refine procedure as in \cite{Dorfler:1996, Verfurth:1996a}. However, compared to adaptive mesh refinement algorithms for the numerical solution for continuous PDEs, the major differences are trifold, and to address those differences, we need new ingredients for the implementation of the adaptive algorithm. 

\begin{itemize}
\item 
The errors $\eta_M$, $\eta_C$ and $\eta_T$ depend on $u_h$ through stress tensors $\sh$ and $\sa$ which are not unique. Therefore, we have to minimize the error estimator with respect to all the admissible stress tensors, and we call this procedure "stress tensor correction". This will be addressed in \S~\ref{sec:Interfacial-stress}.
\item
The truncation error $\eta_T$ is introduced by the truncation of an infinite lattice to a finite domain. If the size of the computation domain is fixed, we shall see the saturation of the numerical error when the degrees of freedom $N$ keep increasing. Therefore, when $\eta_T$ is dominant in the overall error $\eta$, we need to enlarge the computational domain in order to achieve the optimal convergence rate. This will be addressed in \S~\ref{sec:numerics:size}.
\item
The modeling error $\eta_M$ results from the inconsistency of the atomistic model and the continuum model at the interface and finite element edges. In particular, when the interface error is large, we need to enlarge the atomistic domain $\Oma$, and adjust the triangulation in the continuum domain such that the mesh in the continuum region aligns with the micro-triangulation $\Ta$ close to the interface, and the overall triangulation still maintains good quality. This will be addressed in Remark \ref{rem:triage}.
\end{itemize}

\subsubsection{Stress tensor correction}
\label{sec:Interfacial-stress}
By Theorem \ref{thm:residue} and Remark \ref{rem:stressform},  the error estimators $\eta_T$, $\eta_M$, and $\eta_C$ depend on the stress tensors $\sh$ and $\sa$, which are unique up to divergence free tensor fields. Therefore, we need to minimize $\eta(y_h) = \eta_T(y_h)+\eta_M(y_h)+\eta_C(y_h)$ with respect to all the admissible stress tensors. Recall the "stress tensor correction" of the residual estimate \eqref{eq:sharpresidue}, 
\begin{equation}
		 \< \del\Ea(\Ia y_h), v \> \leq  \min_{c_a\in N_1(\T_a)^2, c_h\in N_1(\T_h)^2}\tilde{\eta}(\sa(\Ia y_h)+\nabla c_a\mJ, \sh(y_h)+\nabla c_h\mJ)\|\nabla v\|_{L^2}.
		\label{eq:sharpresidual2}
\end{equation}

In \eqref{eq:sharpresidual2}, we need to solve a nonlinear minimization problem with respect to $c_a$ and $c_h$ which are both defined over whole $\Omega$, the dimension of $c_a$ is $2|\Fs_h|$, and the dimension of $c_h$ is $2|\Fs_\a|$. The cost for the exact stress tensor correction is proportional to solving the original energy minimisation problem.  

Here, we introduce an approximate version of stress tensor correction, which is motivated by the explicit calculation in \cite[Lemma 5.2]{PRE-ac.2dcorners} as well as the analysis of a/c stress tensor in \cite[\S~6.2.3]{Or:2011a}: a "good" a/c stress tensor can be chosen such that it equals to the atomistic stress tensor in the atomistic domain, and equals to the continuum stress tensor for uniform deformation. To be precise, we only need to apply the stress tensor correction to the modelling error $\eta_M$; and in addition, we choose $c_a\equiv 0$, and $c_h(q_f) =0$, where $q_f$ is the midpoint of $f\in \Fh$,  $f\bigcap \L_\i= \emptyset$. Thus the only degrees of freedom to be determined are those $c_h(q_f)$ such that $f\bigcap \L_\i \neq \emptyset$. 

%\textcolor{red}{Notice that, for truncation error \eqref{eqn:etat},  Since we apply the $\eta_T$ on the boundary and $\eta_M$ around the interface, the correction will not conflict.}s

We propose the following algorithm for approximate stress tensor correction: 
\begin{algorithm}
\caption{Approximate stress tensor correction}
\label{alg:etc}
\begin{enumerate}
\item Take $\sa(\Ia y_h)$ and $\sh(y_h)$ as the canonical forms in \eqref{eq:defn_Sa} and \eqref{eq:defn_Sh} respectively. 
\item Denote $q_f$ as the midpoint of $f\in \Fh$. $c_h$ minimizes the following sum
\begin{equation}
	\sum_{T\in\T^\i}|T|\left[\sa(\Ia y_h, T) - \big(\sh(\Ia y_h, T)+\nabla c_h \mJ\big)\right]^2
\end{equation} 
subject to the constraint that $c_h(q_f)=0$, for $f\bigcap \L_\i = \emptyset$.
\item Let $\sh(y_h)  = \sh(y_h) + \nabla c^h\mJ$, compute $\eta_M$, $\eta_T$ and $\eta_C$ with $\sa(\Ia y_h)$ and $\sh(y_h)$.
\end{enumerate}
\end{algorithm}

Instead of minimizing the total error estimator $\eta$ with respect to $c_a$ and $c_h$ as in \eqref{eq:sharpresidual2}, now we only need to minimize the modeling error $\eta_M$ with respect to the degrees of freedom of $\sh$ adjacent to the interface. This dramatically reduced the computational cost of "stress tensor correction". In the implementation, the cost of stress tensor correction is only a small fraction of the total cost, but it greatly improves the accuracy.

\begin{figure} 
	\centering 
	\subfloat[$H_1$ error.]{
		\label{figs:cfdivfree:a}
		\includegraphics[width=6.5cm]{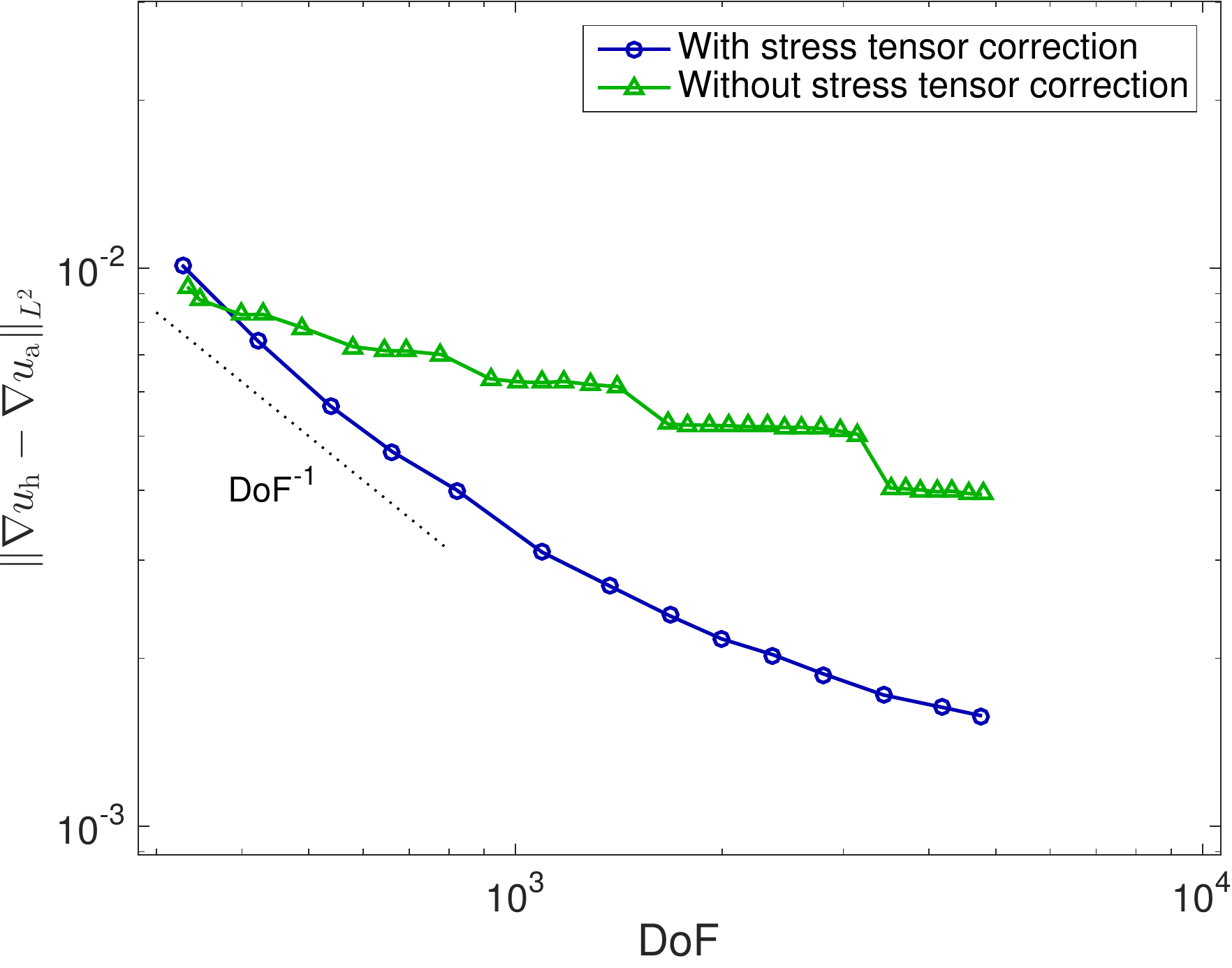}} 
		\hspace{0.5cm} 
	\subfloat[Estimator]{ 
		\label{figs:cfdivfree:b} 
		\includegraphics[width=6.5cm]{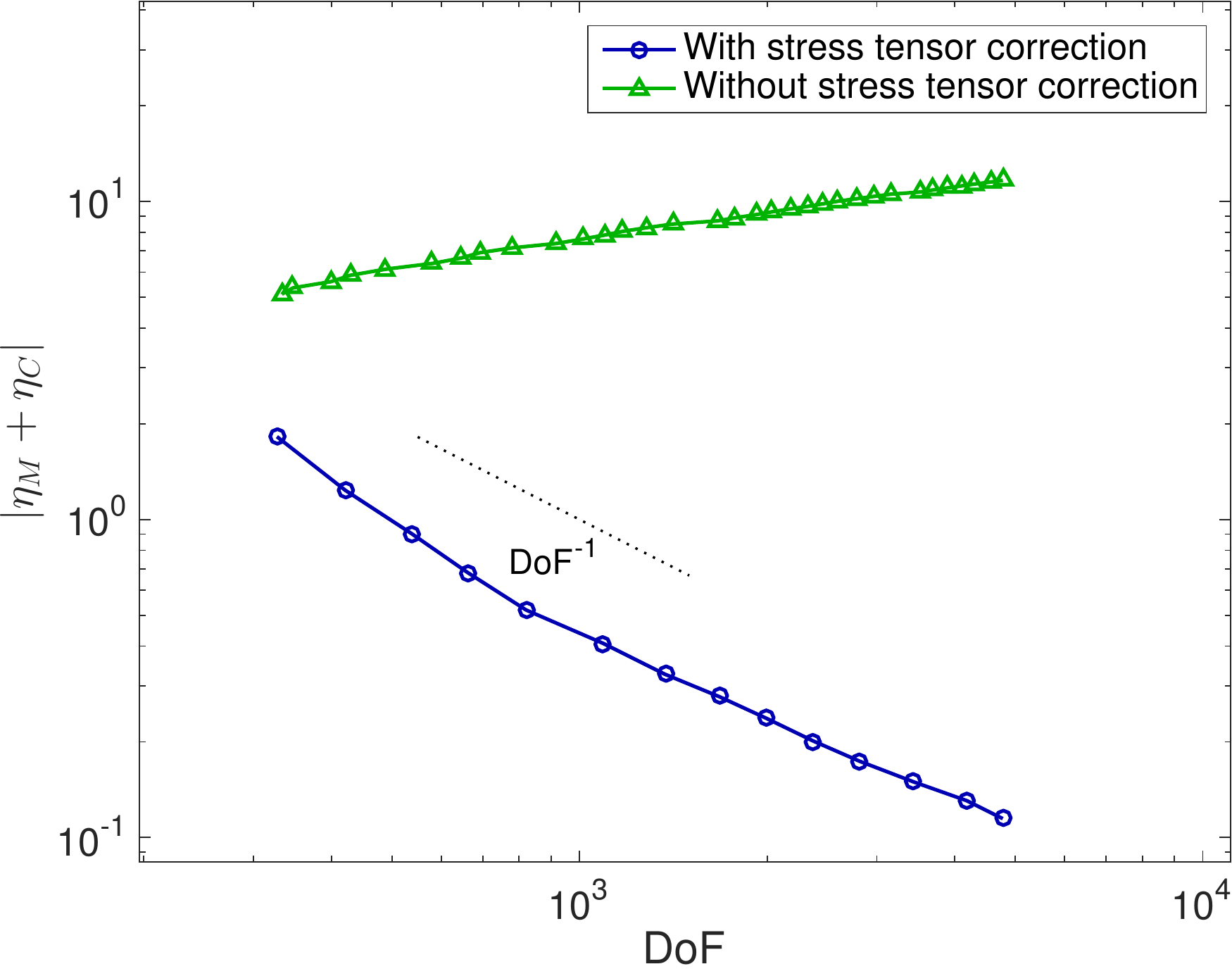}}
		 \caption{Effect of approximate stress tensor correction. Divacancy example, $R_{c}=1000$, take $\tau_1=0.7$ and 
		 $\tau_2=0.2$ in Algorithm \ref{alg:main}. Figure \ref{figs:cfdivfree:a}: $H_1$ error vs. DoF; Figure \ref{figs:cfdivfree:b}: $\eta_M+\eta_C$ vs. DoF.} 
		\label{figs:cfdivfree}
\end{figure}

We numerically demonstrate the effect of the approximate stress tensor correction in Figure \ref{figs:cfdivfree}. We fix the computational domain in this example, therefore we expect the "optimal" error will follow the $N^{-1}$ asymptotics as the degrees of freedom $N$ increase, and get saturated at the level of the truncation error. Figure \ref{figs:cfdivfree:a} shows $H^1$ errors with respect to degrees of freedom $N$. If the stress tensor correction is applied, the error follows the optimal $N^{-1}$ asymptotics before the saturation is reached; if the stress tensor correction is not applied, the error is suboptimal. Figure \ref{figs:cfdivfree:b} shows the error estimator $\eta$ with respect to degrees of freedom $N$. The $N^{-1}$ convergence of $\eta$ is much more significant with correction; without correction $\eta$ may even increase with respect to $N$.

\subsubsection{Local error estimator}
We need to assign global estimators to local elements properly, then mark and subdivide those elements which contribute most to the estimator.     

Recall the definition of $\eta_M$ in \eqref{eqn:etam}, and after taking the stress tensor correction in Algorithm \ref{alg:etc}, we have
\begin{displaymath}
	(\eta_M(y_h))^2 := (C^{\Tr})^2\sum_{T\in\Ta}|T|\big[\sa(\Ia y_h, T) - \sum_{T'\in \Th, T'\bigcap T\neq \emptyset}\frac{|T'\bigcap T|}{|T|}(\sh(y_h, T'))\big]^2.
\end{displaymath}
The contribution is 0 for those $T\in\Ta$ located completely inside an element $T'\in \Th$. As a result, we need only take care of those $T\in\Ta$ and $T'\in\Th$ with $T\bigcap \partial T'\neq \emptyset$. We first define
\begin{displaymath}
\eta_{M}(T, T'):= |T'\bigcap T|\left[\sa(\Ia y_h, T) - \frac{|T'\bigcap T|}{|T|}(\sh(y_h, T'))\right]^2. 
\end{displaymath}
for $T\in\Ta$, then let $\eta_{M}(T') = \sum_{T\in \Ta, T\bigcap T'\neq \emptyset} \eta_M(T, T')$ for $T'\in\Th$. Notice that $(C^{\Tr})^2\sum_{T\in \Th} \eta_M(T)=\eta_M^2$. 
 
Analogously, we can define the local contribution of the truncation error $\eta_T(T')$ for $T'\in \Th$, such that $\sum_{T'\in\Th}\eta_T(T') = \eta_T^2$. Please also refer to Remark \ref{rem:etaT}.

For the coarsening error, recall the definition \eqref{eqn:etac}, 
\begin{displaymath}
\eta_C(y_h) := \sqrt{3}C^{\Tr}C'_{\Th}(\sum_{f\in\Fh} (h_f\sjump)^2)^\frac12,
\end{displaymath}
we define $\eta_{C}(T)$ as follows,
\begin{displaymath}
\eta_{C}(T) = \sqrt{3}C^{\Tr}C'_{\Th}\sum_{f\in \Fh\bigcap T\in\Th}\frac12(h_f\sjump_f)^2.
\end{displaymath}
%  $\frac{|T'_{1}|}{|T'_{1}\cap T'_{2}|}$ and $\frac{|T'_{2}|}{|T'_{1}\cap T'_{2}|}$ of the value computed, respectively. 

For the energy estimator $\mu_E$ from section \S~\ref{sec:error:mainresults:energyerror}, similar to the case of $\eta_M$, we can define the local contributions similarly as $\mu_E(T)$ such that $\sum_{T'\in \Th}\mu_E^2(T') = \mu_E^2$.
\def\M{\mathcal{M}}

Once all the local estimators are assigned, we are ready to define the indicator $\rho_{T}$:
\begin{equation}
\rho_{T} = (C^{\Tr})^2\frac{\eta_{M}(T)}{\eta_M}+(C^{\Tr})^2\frac{\eta_{T}(T)}{\eta_T}+(\sqrt{3}C^{\Tr}C'_{\Th})^2\frac{\eta_{C}(T)}{\eta_C}.
\label{eq:localestimator1}
\end{equation}
Notice that the sum of local estimators is equal to the global estimator.

Meanwhile, for the energy based estimate, we have,
% (\lz{Mingjie, please check})
\begin{equation}
\rho^E_{T} = C^E(C^{\Tr})^2\left(\eta_{M}(T)^2+\eta_T(T)^2\right)+C^E(\sqrt{3}C^{\Tr}C'_{\Th})^2\left(\eta_{C}(T)\right)^{2} + |\mu_{E}(T)|
\label{eq:localestimator2}
\end{equation}

The constants $C^{\Tr}$, $C^E$, $C'_{\Th}$ in \eqref{eq:localestimator1} and \eqref{eq:localestimator2} are not known a priori, instead, we use their empirical estimates in the implementation.

Algorithm \ref{alg:main} is the main algorithm for the adaptive mesh refinement, and D\"{o}rfler adaptive strategy \cite{Dorfler:1996} is used in the algorithm.

\begin{algorithm}
\caption{A posteriori mesh refinement}
\label{alg:main}
\begin{enumerate}
	\item[Step 0] Prescible $\Omega_R$, $\Th$, $N_{\max}$, $\rho_{\rm tol}$, $\tau_1$ and $\tau_2$.
	\item[Step 1] \textit{Solve:} Solve the a/c solution $y_h$ of \eqref{eq:min_ac} on the current mesh $\Th$. 
	\item[Step 2] \textit{Estimate}: Carry out the stress tensor tensor correction step in Algorithm \ref{alg:etc}, and compute the error indicator $\rho_{T}$ for each $T\in\Th$. For fixed $R$, we do not need to include the contribution from truncation error $\eta_T$ in $\rho_T$. Set $\rho_{T} = 0$ for $T\in\Ta\bigcap \Th$. Compute the degrees of freedom $N$ and total error $\rho = \sum_{T}\rho_T$. Stop if $N>N_{\max}$ or $\rho < \rho_{\rm tol}$.
	\item[Step 3] Mark: 
	\begin{enumerate}
	\item[Step 3.1]:  Choose a minimal subset $\M\subset \Th$ such that
	\begin{displaymath}
		\sum_{T\in\M}\rho_{T}\geq\frac{1}{2}\sum_{T\in\Th}\rho_{T}.
	\end{displaymath}	 
	\item[Step 3.2]: Find the interface elements $\M_\i:=\{T\in\M: T\bigcap \L_\i\neq \emptyset\}$. Check if
	\begin{equation}
		\sum_{T\in \M_\i}\rho_{T}\geq \tau_1\sum_{T\in\M}\rho_{T}.
		\label{eq:interface1}
	\end{equation}
	where tolerance $0<\tau_1<1$.	If true, let $\M = \M\setminus \M_\i$.
	\end{enumerate}
	\item[Step 4] \textit{Refine:} If \eqref{eq:interface1} is true, expand interface $\L_\i$ outward by one layer. Then, bisect all elements $T\in \M$. Stop if $\frac{\eta_T}{\eta_M + \eta_C}\geq \tau_2$, otherwise, go to Step 1. 
\end{enumerate}
\end{algorithm}

\begin{remark}
For the calculation with fixed computational domain, the numerical error will saturate at the level of truncation error. The stoping criteria can be modified as: 

Step 2: ... Compute the convergence rate $\beta$ of the estimated total error $\rho$ with respect to the degrees of freedom $N$. Stop if $ \beta \leq \tau_2$.
\label{rem:stop}
\end{remark}

\begin{remark}
It is possible to use different mark strategies, for example, 

\textit{Step 3.1}: Choose a minimal subset $\M$, s.t. 
	\begin{displaymath}
		\rho_{T}\geq\textrm{mean}(\rho), \quad\forall T\in\M.
	\end{displaymath}

\def\dist{\rm dist}

\textit{Step 3.2} We can find the interface elements which are within $k$ layers of atomistic distance, $\M^k_\i:=\{T\in\M\bigcap \Th^\c: \dist(T, \Li)\leq k \}$. Choose $K\geq 1$, find the first $k\leq K$ such that 
	\begin{equation}
		\sum_{T\in \M^k_\i}\rho_{T}\geq \tau_1\sum_{T\in\M}\rho_{T},
		%\label{eq:interface2}
	\end{equation}
	with tolerance $0<\tau_1<1$.	If such a $k$ can be found, let $\M = \M\setminus \M^k_\i$. Then in step 3, expand interface $\L_\i$ outward by $k$ layers. 
\end{remark}

\begin{remark}
After pushing the interface outward in Step 4, we have to 'remove' those triangles in the continuum mesh which overlap with the new atomistic region. It will generate a gap between the atomistic region and the continuum region. We need to triangulate this gap, and adjust the positions of the nodes to improve the quality of the interfacial triangles. In our implementation, we adapted the Matlab package \textit{EasyMesh}, a two-dimentional quality mesh generator to carry out this task \cite{easymesh}.
\label{rem:triage}
\end{remark}

\begin{figure} 
	\centering 
	\subfloat{
		\label{figs:easymesh:a}
		\includegraphics[width=5cm]{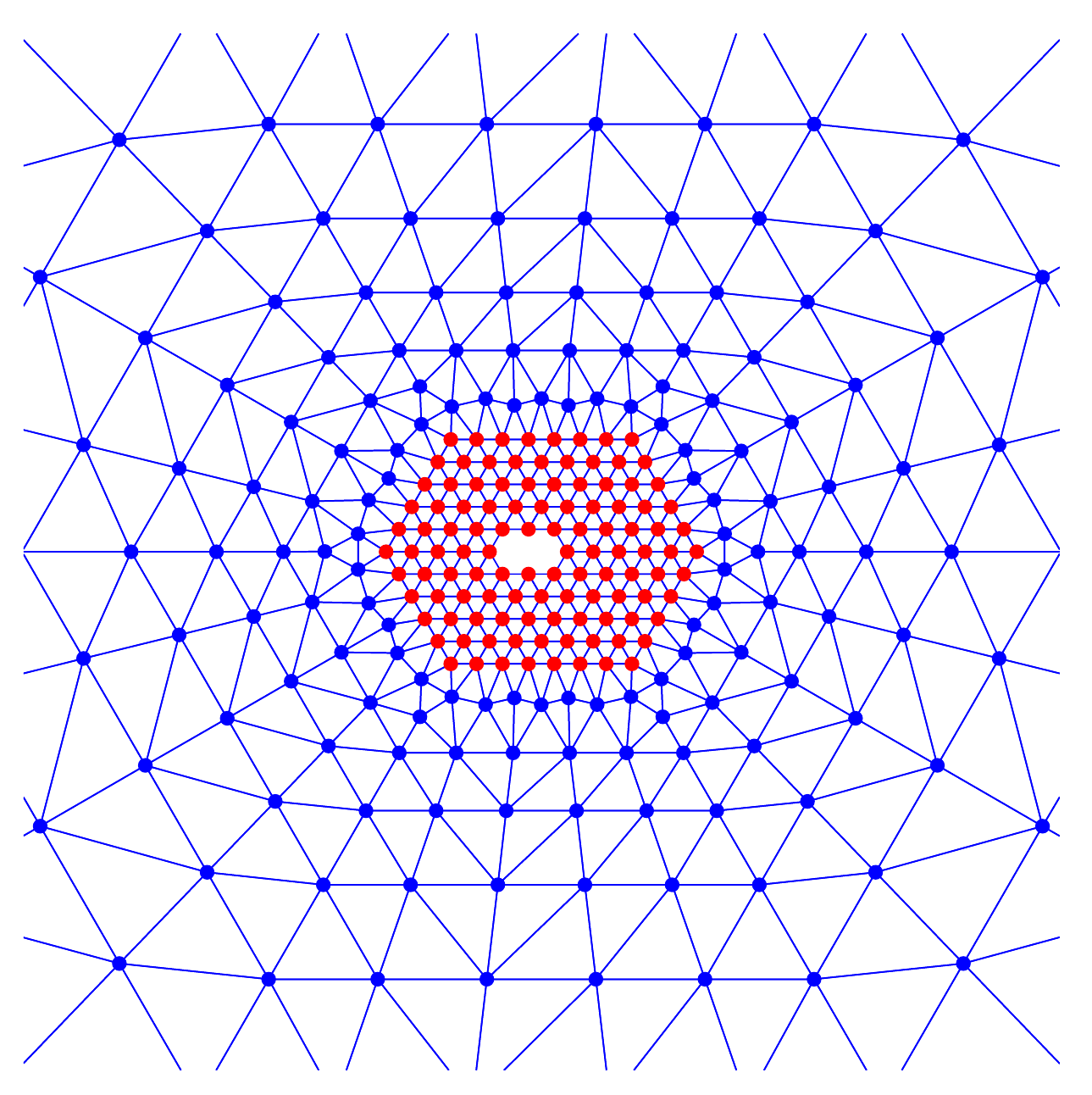}} 
%		\hspace{0.2cm}
	\subfloat{
		\label{figs:easymesh:b}
		\includegraphics[width=5cm]{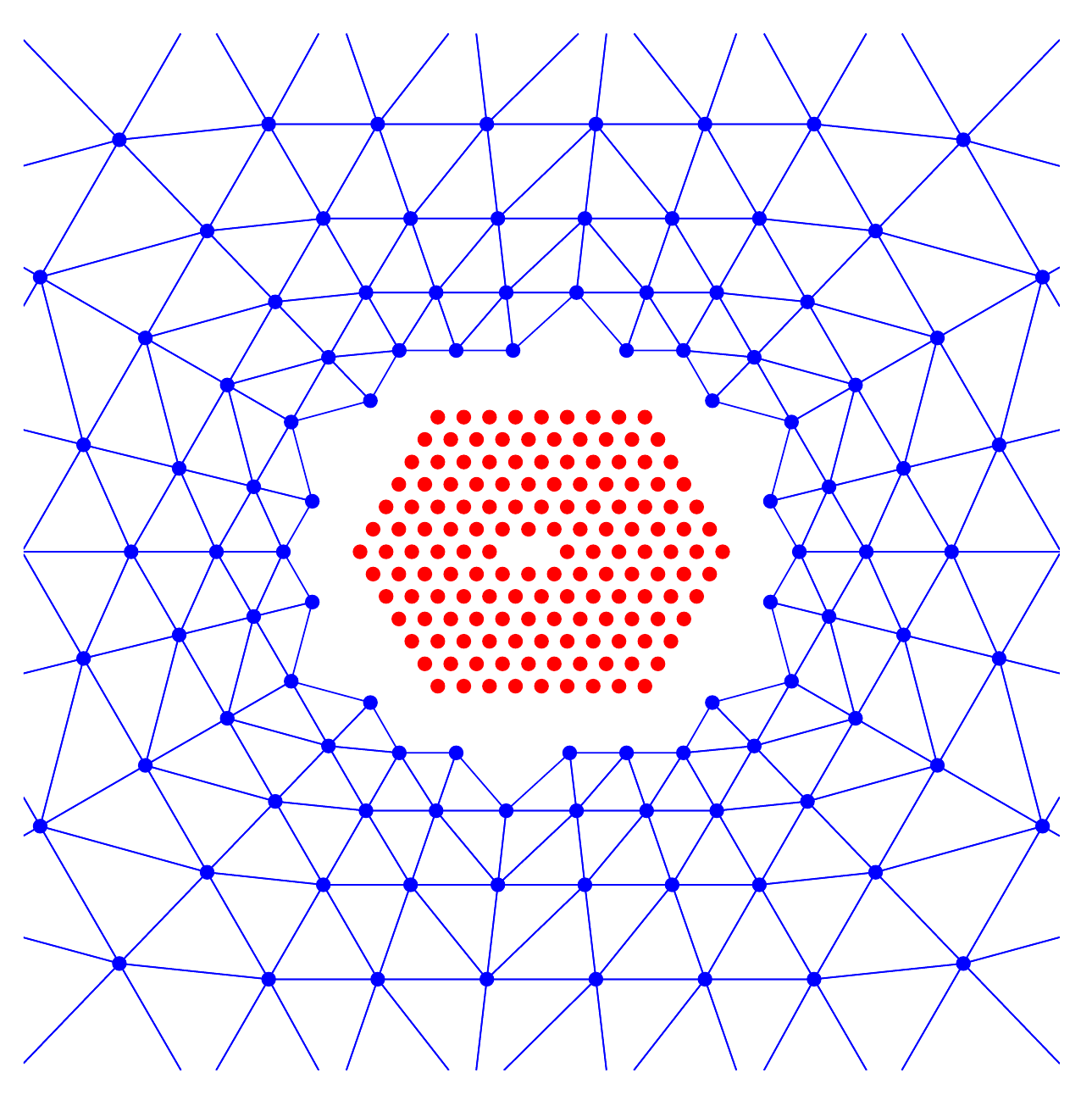}} \\
%		\hspace{0.2cm} 
	\subfloat{ 
		\label{figs:easymesh:c} 
		\includegraphics[width=5cm]{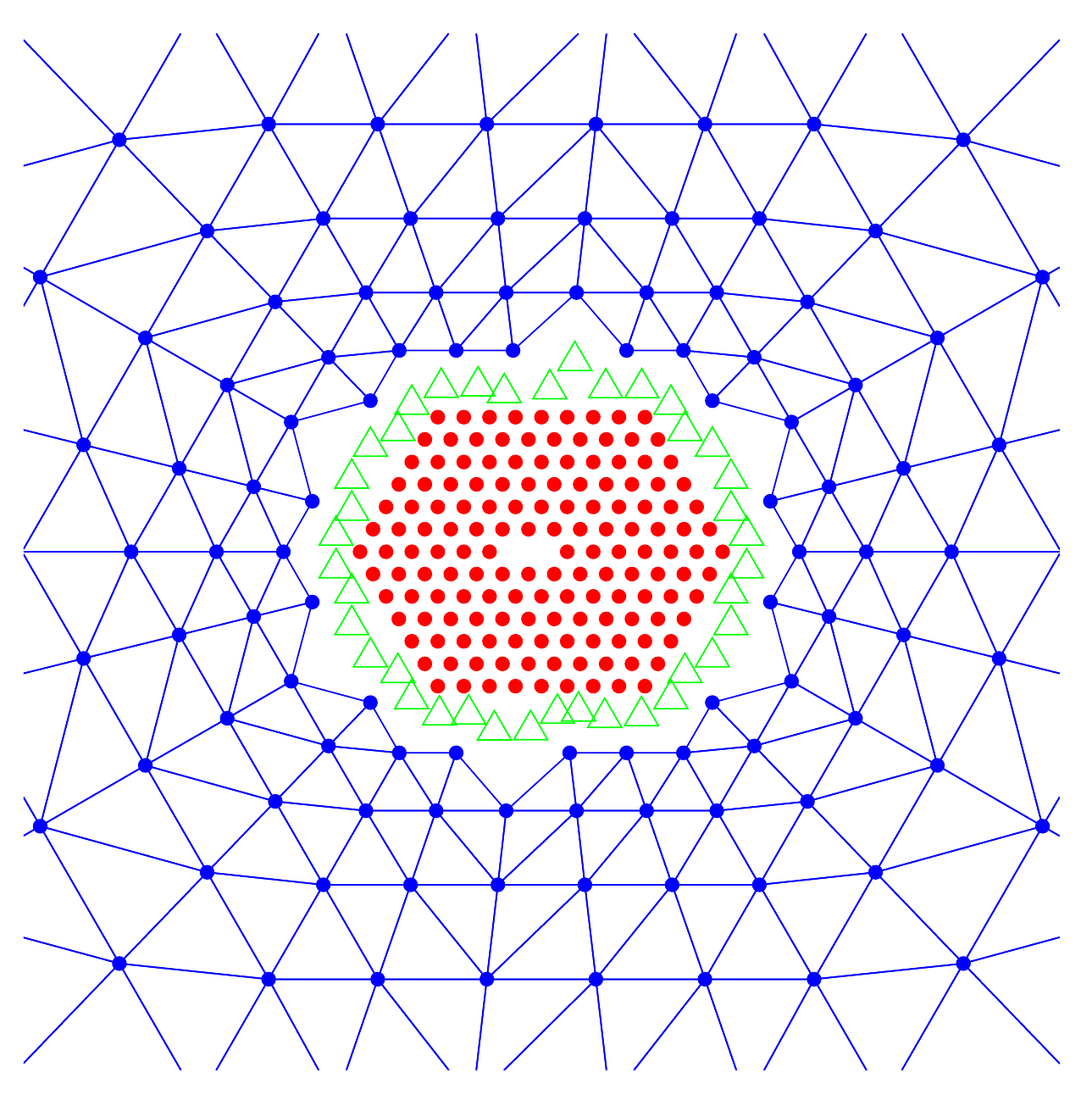}}
%		\hspace{0.2cm}
	\subfloat{
		\label{figs:easymesh:d}
		\includegraphics[width=5cm]{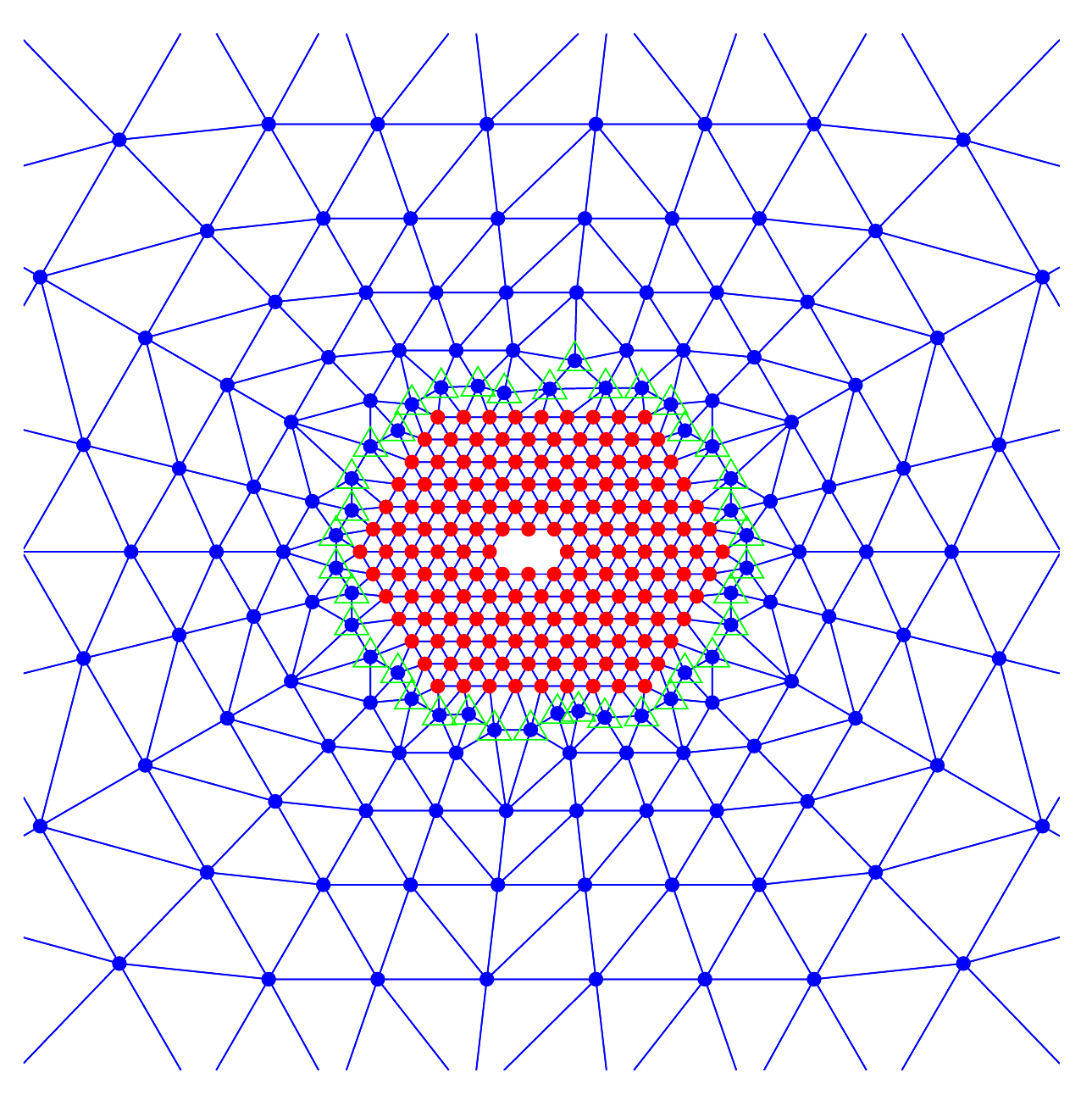}}  
		 \caption{Snapshots of the expanding interface in Step 4 of Algorithm \ref{alg:main}. (Top-left) initial mesh with $R_a=5$; (Top-right) mesh with $R_a=6$: after removing the neighboring continuum nodes close to the interface, move the interface outward by 1 layer ; (Bottom-left) generating new continuum nodes (marked with green triangles) and adjusting their positions to maintain the quality of mesh; (Bottom-right) final triangulations.} 
		\label{figs:easymesh}
\end{figure}

\subsection{Model Problem}
\label{sec:numerics:problem}
%
% Throughout this section, we fix the computational domain size $N=100$, and 
Recall the EAM potential defined in \eqref{eq:eam_potential}. Let
\begin{displaymath}
\phi(r)=\exp(-2a(r-1))-2\exp(-a(r-1)),\quad \psi(r)=\exp(-br)
\end{displaymath}
\begin{displaymath}
F(\tilde{\rho})=C\left[(\tilde{\rho}-\tilde{\rho_{0}})^{2}+
(\tilde{\rho}-\tilde{\rho_{0}})^{4}\right]
\end{displaymath}
with parameters $a=4, b=3, c=10$ and $\tilde{\rho_{0}}=6\exp(0.9b)$, which is the same as the numerical experiments in the a priori analysis paper \cite{COLZ2013}.

\begin{figure}[htb]
\begin{center}
	\includegraphics[scale=0.4]{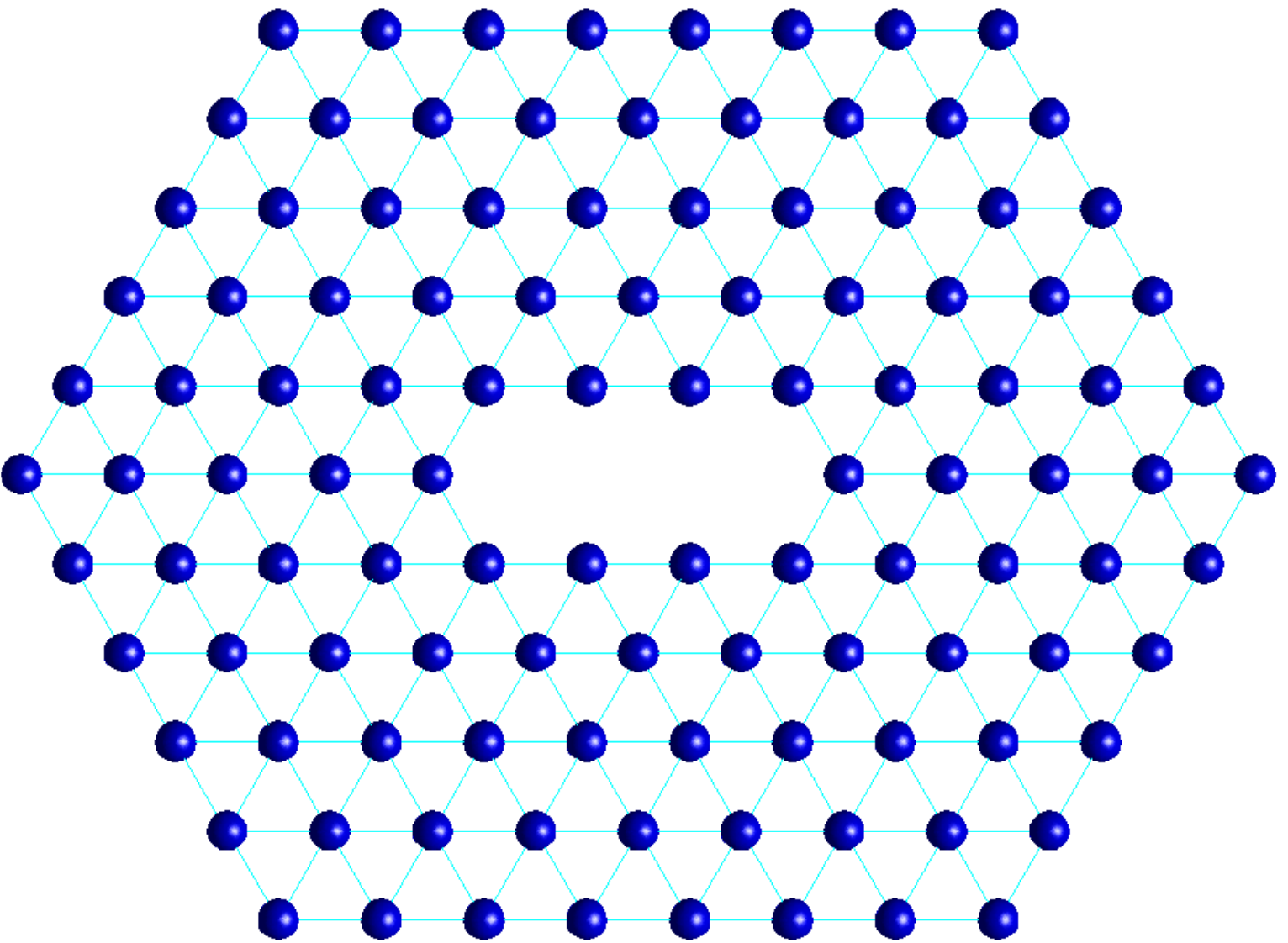}
	\caption{Illustration of the atomistic lattice $\L$ with 3 vacancies surrounded by 5 atomistic layers.}
	\label{figs:3v5l-config}
\end{center}
\end{figure}

To generate a defect, we remove $k$ atoms from $\Lhom$,
\begin{align*}
\L_{k}^{\rm def}:=\{-(k/2)e_{1}, \ldots, (k/2-1)e_{1})\},     & \qquad{\rm if }\quad k \quad\textrm{is even},\\
\L_{k}^{\rm def}:=\{-(k-1)/2e_{1}, \ldots, (k-1)/2e_{1})\}, & \qquad{\rm if }\quad k \quad\textrm{is odd},
\end{align*}
and $\L = \Lhom\setminus \L_{k}^{\rm def}$. See Figure \ref{figs:3v5l-config} for an illustration.

For $\ell\in\L$, consider the nearest neighbour interaction, $\Nhd_{\ell} := \{ \ell' \in \L \sep 0<|\ell'-\ell| \leq 1 \}$, and interaction range $\Rg_\ell := \{\ell'-\ell \sep \ell'\in \Nhd_\ell\} \subseteq \{a_j, j=1,\dots, 6\}$. The defect core $\Ddef$ can be defined by $\Ddef = \{x: \rm{dist} (x,  \L_{k}^{\rm def})\leq 1 \}$, $\L \bigcap \Ddef$ is the first layer of atoms around $\L_{k}^{\rm def}$.

% fig 5 layers around 3 vacancies.

\subsection{Di-vacancy Example}
\label{sec:numerics:di-vacancy}

In this section, we numerically justify the performance of the proposed adaptive mesh refinement algorithm.
We take the same di-vacancy example in \cite{COLZ2013}, namely, setting $k=2$ for $\L_{k}^{\rm def}$. We apply isotropic stretch $\mathrm{S}$ and shear $\gamma_{II}$  by setting
\begin{displaymath}
{\sf B}=\left(
	\begin{array}{cc}
		1+\mathrm{S} & \gamma_{II} \\
		0            & 1+\mathrm{S}
	\end{array}	 \right)
	\cdot{\sf {F_{0}}}
\end{displaymath}
where ${\sf F_{0}} \propto \mathrm{I}$ minimizing the Cauchy-Born energy density $\mathrm{W}$, $\mathrm{S}=\gamma_{II}=0.03$.
In our numerical experiments, the reference solution denoted as $u_r$ is solved by GRAC method with a sufficient large mesh where $R_\a = 93$ and $R = 17298$. 

%For simplicity, we denote $\epsilon^{\rm H1}$ and $\epsilon^{\rm E}$ as the actual H1 error $\|\nabla u_h - \nabla u_r \|_{L^2}$ and energy  difference ($\|\Eh - \E^r \|_{L^2}$)) with $u_h$ solved by residual estimator driven algorithm, $\epsilon^{\rm H1}_\E$ and $\epsilon^{\rm E}_\E$ the H1 error and energy difference with solutions solved by energy estimate driven algorithm. 

\subsubsection{Fixed computation domain}
\begin{figure}[htp]
\begin{center}
	\includegraphics[scale=0.45]{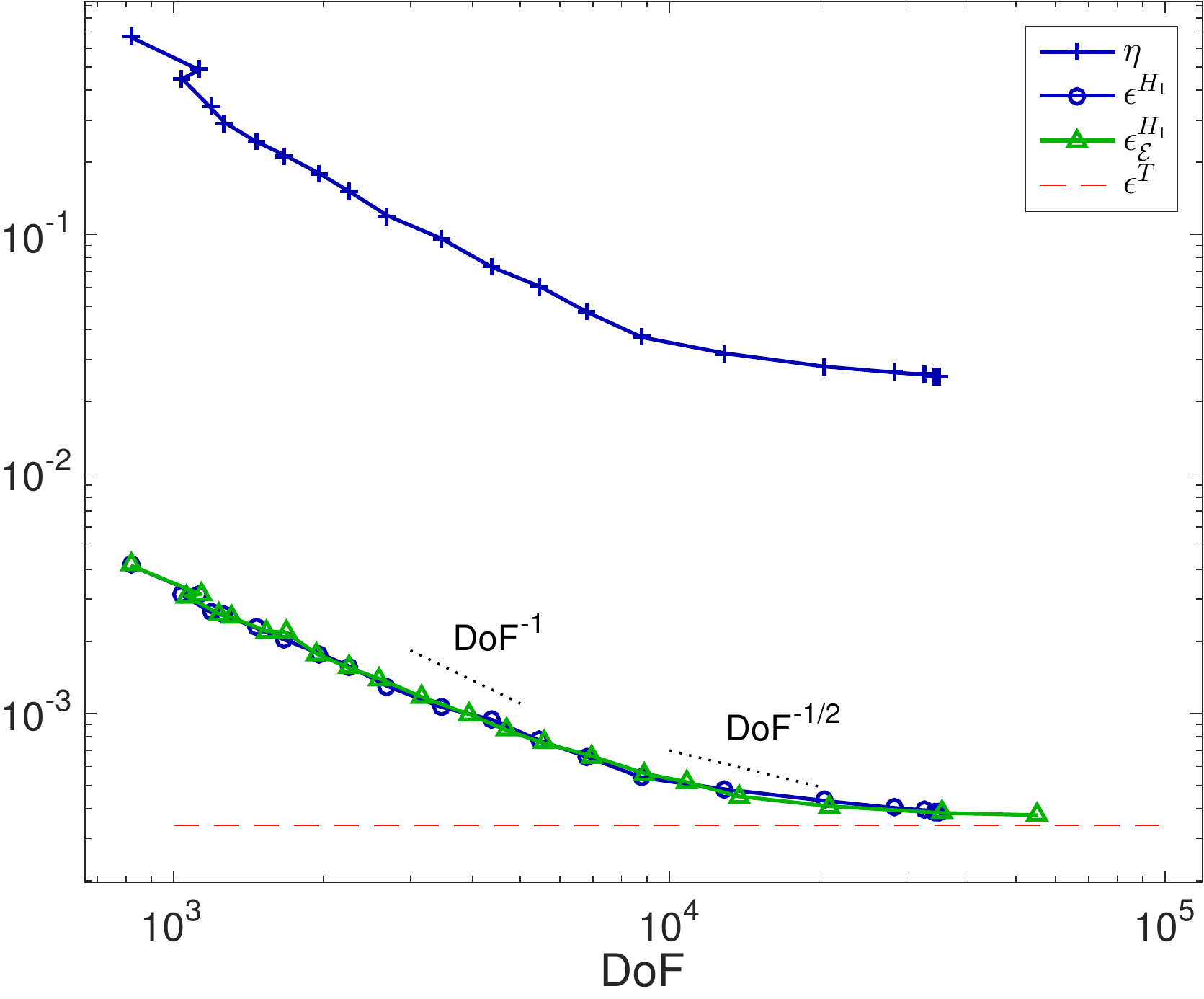}
	\caption{Numerical results by Algorithm \ref{alg:main} and Remark \ref{rem:stop} with $R=1000$, $\tau_1=0.7$, $\tau_2=0.2$. we denote $\epsilon^{{\rm H}_1}$ as the actual H1 error $\|\nabla u_h - \nabla u_r \|_{L^2}$ with $u_h$ solved by residual estimator driven algorithm, $\epsilon^{{\rm H}_1}_{\E}$ as the H1 error with solutions solved by energy estimate driven algorithm, $\epsilon^{T}$ the actual residual  truncation error.}
	\label{figs:H1-1000}
\end{center}
\end{figure}

\begin{figure}[htp]
\begin{center}
	\includegraphics[scale=0.45]{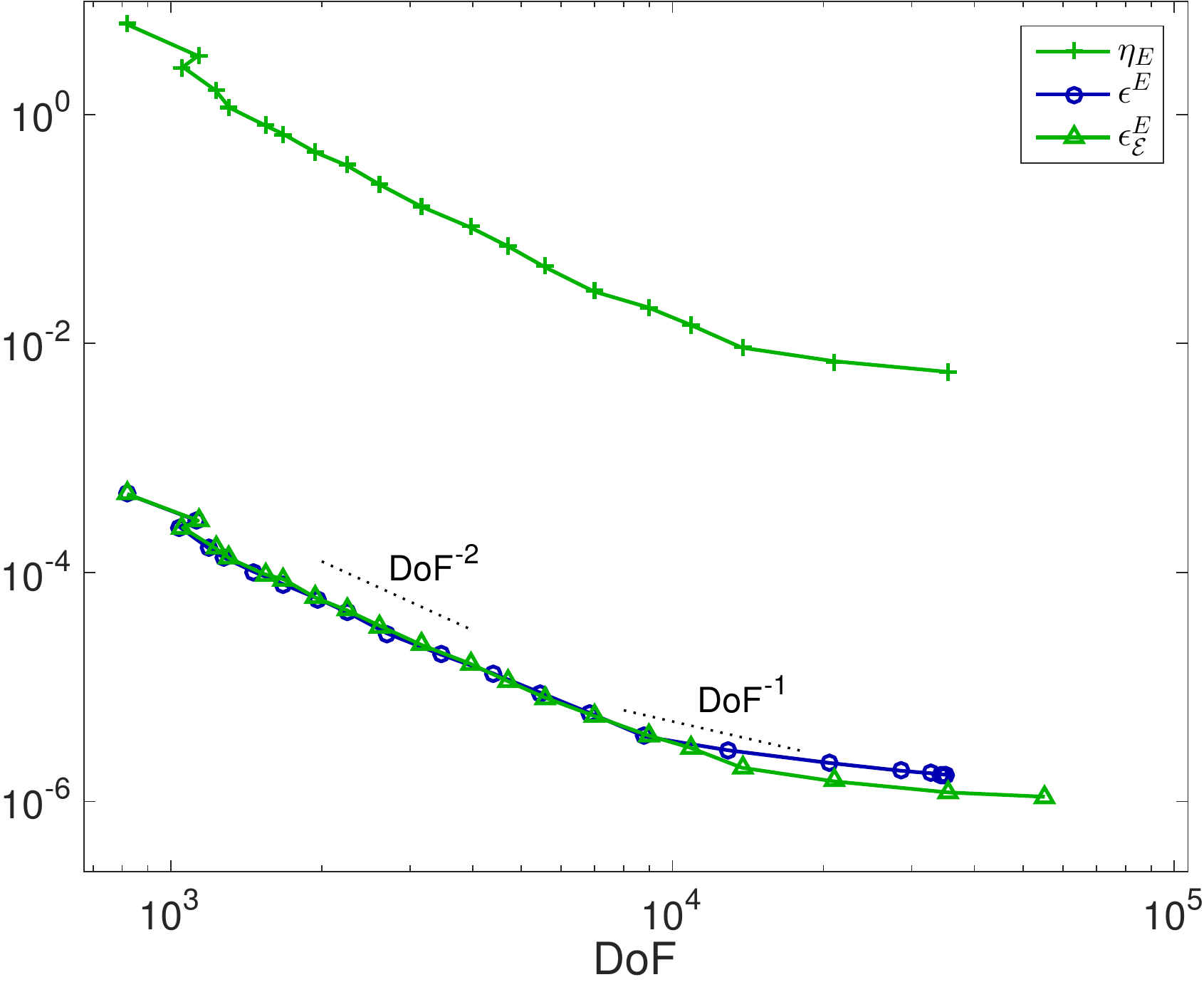}
	\caption{Numerical results by Algorithm \ref{alg:main} and Remark \ref{rem:stop} with $R=1000$, $\tau_1=0.7$, $\tau_2=0.2$. we denote $\epsilon^{\rm E}$ as the actual energy difference $\|\Eh - \E^r \|_{L^2}$ with $u_h$ solved by residual estimator driven algorithm, $\epsilon^{\rm E}_\E$ as the energy difference with solutions solved by energy estimate driven algorithm, $\epsilon^{T}_{\E}$ the actual energy  truncation error.}
	\label{figs:En-1000}
\end{center}
\end{figure}

In this subsection, we fix $R=1000$. The numerical results are shown in \fig{figs:H1-1000} and \fig{figs:En-1000}. The red dashed lines in both figures denote the truncation errors $\eta_T$ and $\eta^2_T$ respectively. The figures show that when $N$ is small, the modelling error and coarsening error dominates, our results coincide with the optimal a priori convergence rate ($N^{-1}$ for $H^1$ norm and $N^{-2}$ for energy, respectively). When $N$ increases, the truncation error becomes dominant, which results in a suboptimal convergence rate and finally saturates the overall error. These results indicate that for a fixed computational domain, we can only achieve optimal convergence rate up to a certain critical degree of freedom. A possible cure is to enlarge the computational domain in order to balance the truncation error with the modeling and coarsening errors, which motivates the next numerical experiments.

% fig H1-1000-Tr
% fig En-1000-Tr

\subsubsection{Adaptive algorithm with automatic control on domain size}
\label{sec:numerics:size}

With the estimator $\eta_T$ for the truncation error, we can modify the Algorithm \ref{alg:main} to automatically enlarge the computational domain if the truncation error is dominant in the total error $\rho$. 

\begin{algorithm}
\caption{A posteriori mesh refinement with size control.}
\label{alg:size}
\begin{enumerate}
	\item[Step 0] Prescible $\Omega_{R_0}$, $\Th$, $N_{\max}$, $\rho_{\rm tol}$, $\tau_1$, $\tau_{3}$ and $R_{\max}$.
	\item[Step 1] \textit{Solve:} Solve the a/c solution $u_{h, R}$ of \eqref{eq:min_ac} on the current mesh $\T_{h, R}$.  
	\item[Step 2] \textit{Estimate}: carry out the stress tensor tensor correction step in Algorithm \ref{alg:etc}, and compute the error indicator $\rho_{T}$ for each $T\in\Th$, including the contribution from truncation error $\eta_T$. Set $\rho_{T} = 0$ for $T\in\Ta\bigcap \Th$. Compute the degrees of freedom $N$, error estimator $\rho_T$ and $\rho = \sum_{T}\rho_T$. Stop if $N>N_{\max}$ or $\rho < \rho_{\rm tol}$ or $R>R_{\max}$.
	\item[Step 3] Mark: 
	\begin{enumerate}
	\item[Step 3.1]:  Choose a minimal subset $\M\subset \Th$ such that
	\begin{displaymath}
		\sum_{T\in\M}\rho_{T}\geq\frac{1}{2}\sum_{T\in\Th}\rho_{T}.
	\end{displaymath}	 
	\item[Step 3.2]: We can find the interface elements which are within $k$ layers of atomistic distance, $\M^k_\i:=\{T\in\M\bigcap \Th^\c: \rm{list}(T, \Li)\leq k \}$. Choose $K\geq 1$, find the first $k\leq K$ such that 
	\begin{equation}
		\sum_{T\in \M^k_\i}\rho_{T}\geq \tau_1\sum_{T\in\M}\rho_{T},
		\label{eq:interface2}
	\end{equation}
	with tolerance $0<\tau_1<1$.	If such a $k$ can be found, let $\M = \M\setminus \M^k_\i$. Then in step 3, expand interface $\L_\i$ outward by $k$ layers. 
	\end{enumerate}
	\item[Step 4] \textit{Refine:} If \eqref{eq:interface2} is true, expand interface $\L_\i$ outward by one layer. If $\eta_{T}\geq \tau_3 \rho$, enlarge the computational domain (details in Remark \ref{rem:enlarge}) . Bisect all elements $T\in \M$. Go to Step 1.
\end{enumerate}
\end{algorithm}

\begin{remark}
In our current implementation, we first generate an initial graded triangulation on $\Omega_{R_{\max}}$ in a way that it contains the triangulation of a sequence of domains $\Omega_{R_k}$ such that $R_0<R_1<\cdots< R_{\max}$. Therefore, when we need to enlarge the computational domain in Step 4 of the above algorithm, we simply combine the triangulation for the current domain $\Omega_{R_k}$ and the initial triangulation of $\Omega_{R_{k+1}}\setminus \Omega_{R_k}$ to generate the triangulation for $\Omega_{R_{k+1}}$. 
\label{rem:enlarge}
\end{remark}

\begin{figure}[!htp]
\begin{center}
	\includegraphics[scale=0.45]{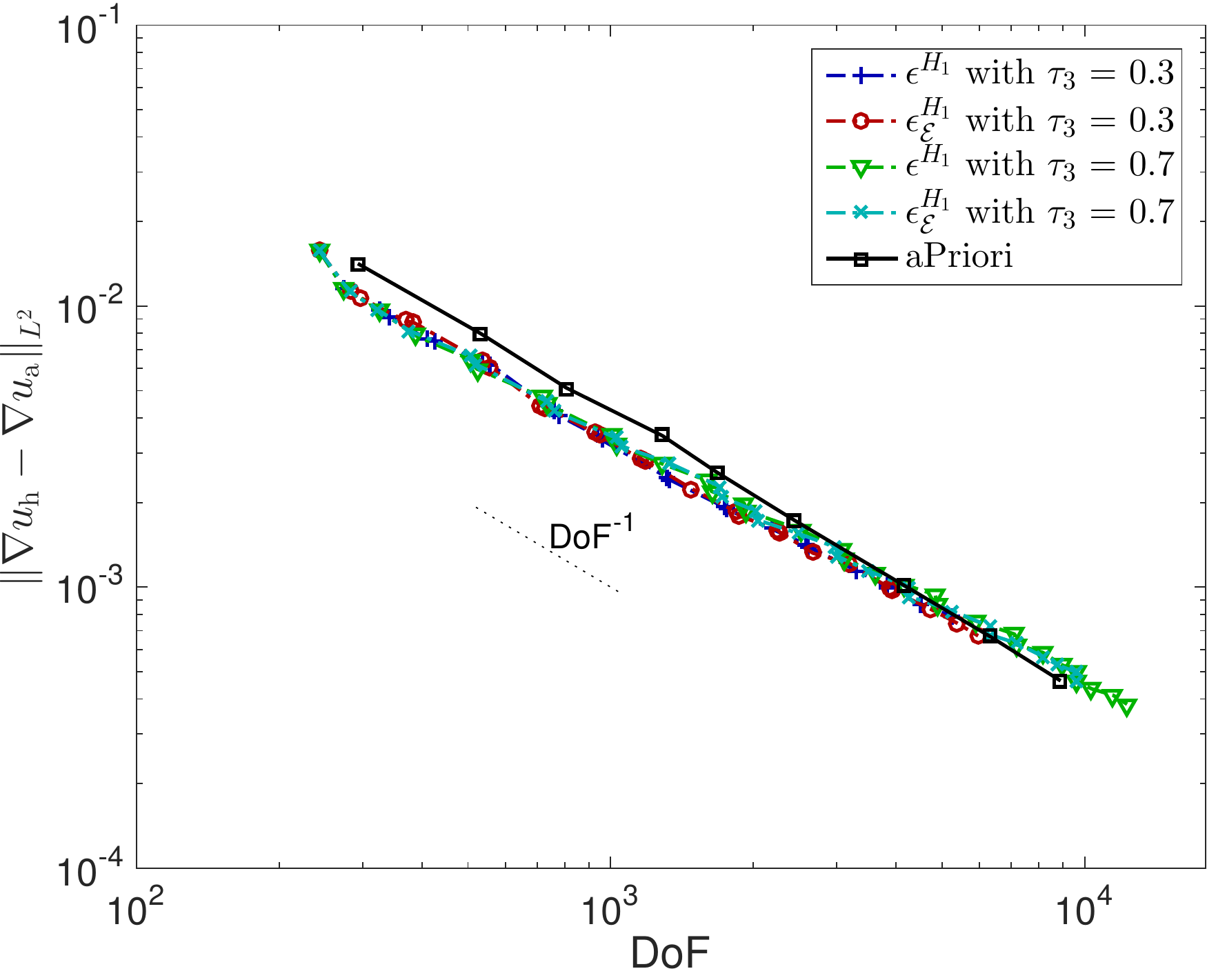}
	\caption{Numerical results by Algorithm \ref{alg:size} and Remark \ref{rem:enlarge}: $H^1$ error vs. Degree of Freedom with $\tau_3 = 0.3$ and $\tau_3 = 0.7$ or both residual estimate driven and energy estimate driven algorithms. The aPriori curve shows the corresponding a priori convergence. }
	\label{fig:changeR-H1}
\end{center}
\end{figure}

\begin{figure}[!htp]
\begin{center}
	\includegraphics[scale=0.45]{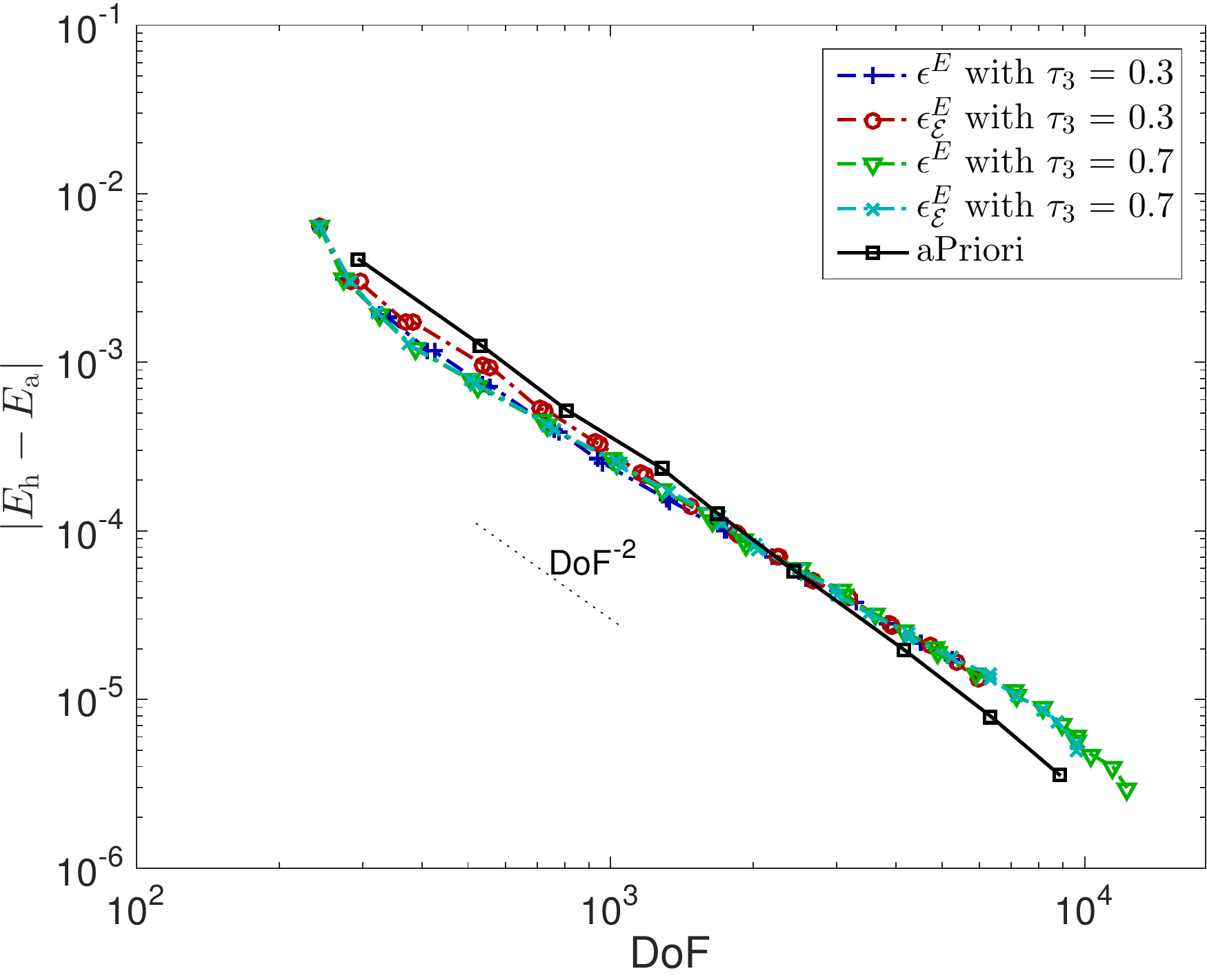}
	\caption{Numerical results by Algorithm \ref{alg:size} and Remark \ref{rem:enlarge}: Energy difference vs. Degree of Freedom with $\tau_3 = 0.3$ and $\tau_3 = 0.7$ for both residual estimate driven and energy estimate driven algorithms. The aPriori curve shows the corresponding a priori convergence.}
	\label{fig:changeR-En}
\end{center}
\end{figure}

From the numerical results in Figures \ref{fig:changeR-H1} - \ref{fig:changeR-En}, we can see that with Algorithm \ref{alg:size}, it is possible to change the domain size automatically, and maintain the optimal convergence rate without the error saturation phenomenon we observed for fixed size computations. The parameter $\tau_3$ can be used to tune the balance between truncation error and other error contributions. With a smaller $\tau_3$, the algorithm tends to enlarge the domain more frequently, while with a larger $\tau_3$, the algorithm tends to push outward the atomistic region and refine the coarse mesh more frequently. In the numerical results, we test two values $\tau_3 = 0.3$ and $\tau_3 = 0.7$. Although there are some small differences, the overall convergence behaviour looks similar and are comparable to the a priori results. 

%\begin{figure}[!htp]
%\begin{center}
%	\includegraphics[scale=0.45]{./ED-H1.pdf}
%	\caption{Numerical results by Algorithm \ref{alg:main} and Remark \ref{rem:enlarge}: $H^1$ error vs. $N$, using a posteriori estimator in energy.}
%	\label{fig:changeR-H1-En}
%\end{center}
%\end{figure}
%
%\begin{figure}[!htp]
%\begin{center}
%	\includegraphics[scale=0.45]{./ED-En.pdf}
%	\caption{Numerical results by Algorithm \ref{alg:main} and Remark \ref{rem:enlarge}: Energy error vs. $N$, using a posteriori estimator in energy.}
%	\label{fig:changeR-En-En}
%\end{center}
%\end{figure} 

\clearpage
\newpage

\section{Conclusion}
\label{sec:conclusion}
In this paper, we derive rigorous a posteriori error estimates for a class of consistent (ghost force free) atomistic/continuum coupling schemes. Numerical results for the corresponding adaptive algorithms are comparable to optimal a priori analysis. This opens an avenue for further mathematical analysis and algorithmic developments for longer range interactions, higher dimensional problems, and general atomistic/continuum coupling algorithms. 

For general short range interactions longer than the nearest neighbour, the stress tensor can be defined using the localization formula and quasi-interplant as in the a priori analysis \cite{Or:2011a, OrShSu:2012, OrtnerTheil2012}. The residual estimate can be carried out analogously as in this paper. However, such a stress tensor is not anymore piecewise constant, and may require complicated geometric operations to evaluate. Therefore, the numerical implementation is difficult and we are currently pursuing an alternative approach to define piecewise constant stress tensor field for general short range interactions. 

The extension to the case of the straight screw dislocation in 2D and point defect case in 3D is straightforward. More practical problems, for example, the study of dislocation nucleation and dislocation interaction by a/c coupling methods has attracted considerable attention from the early stage of a/c coupling methods \cite{Tadmor:1999, Phillips:1999}. The difficulty is to deal with boundary condition and complicated geometry changes of the interface.

For general atomistic/continuum coupling schemes, such as BQCE, BQCF and BGFC, the a priori analysis in \cite{MiLu:2011, LiOrShVK:2014, OrZh:2016} provide a general analytical framework and the stress tensor based formulation plays a key role in the analysis. Therefore, the a posteriori analysis for those coupling schemes can inherit this analytical framework and the stress tensor formulation. The stress tensor correction method and other techniques developed in this paper will be essential for the efficient implementation of the corresponding adaptive algorithms.

\subsection*{Acknowledgement}

The authors thank Christoph Ortner and Huajie Chen for the stimulating discussions on the adaptive computation of material defects. The authors also thank the referees for their insightful comments. Their feedback has helped clarify various aspects of our work. 

%%%%%%%%%%%%%%%%%%%%%%%%%%%%%%%%%%%%%%%%%%%%%%%%%%%%
%Appendix%
%%%%%%%%%%%%%%%%%%%%%%%%%%%%%%%%%%%%%%%%%%%%%%%%%%%%
\appendix
% \section{Appendix}
%
% \label{sec:appendix}
%
\section{Extension to the vacancies}
\label{sec:appendix:extension}

\def\PhiB{\Phi_{\bbB}}

We need to extend $v$ from $\L$ to $\Lhom$ which includes the vacancy sites. We first define the extension operator $E$ on $\Us$ by 
\begin{equation}
Eu := \argmin_{v \in \Us, v = u \text{ on } \L} \PhiB(v):=
\argmin \sum_{b \in \bbB} |{\rho}_b\cdot D_b v|^2, \quad
\forall u \in \Us,
\end{equation}
where $\bbB$ defined in \eqref{def: definition of bonds including vacancies}
is the set of all nearest-neighbour interaction bonds in $\Lhom$. 
Notice that for $v\in \Us$, $\| \nabla v\|_{L^2}$ can be properly and uniquely defined by 
$\| \nabla Ev\|_{L^2}$.  

It is known from \cite[Proposition 4.1]{OrtnerShapeev:2011} that $\PhiB(v)$ is equivalent to $\|\nabla v\|_{L^2}$ such that, 
\begin{equation}
	\frac{3}{4}\|\nabla v \|_{L^2}^2\leq \PhiB(v) \leq \frac{9}{4}\|\nabla v \|_{L^2}^2
	\label{eq:phibequiv}
\end{equation}

Since $\mA^{-1}E\mA v = Ev$ on $\L$, by definition of $Ev$, we have $\PhiB(\mA^{-1}E\mA v) \geq \PhiB(Ev)$. Combining with the inequality $\|\mG\mH\|_F\leq \|\mG\|_F\|\mH\|_F$ for the matrix Frobenius norm and \eqref{eq:phibequiv}, it holds that,
\begin{align}
	\|\nabla E v\|_{L^2}^2 & \leq \frac{4}{3} \PhiB(E v)\nonumber\\
			       & \leq \frac{4}{3}\PhiB(\mA^{-1}E\mA v)\nonumber\\
			       & \leq 3 \|\nabla\mA^{-1}E\mA v\|_{L^2}^2\nonumber\\
			       & \leq 3 \|\mA^{-1}\|_F\|\nabla E\mA v\|_{L^2}^2. \label{eq:extensionineq}
\end{align}

\bibliographystyle{plain}
\bibliography{qc-1.bib}
\end{document}